\documentclass[11pt]{article} 
\usepackage{authblk}
\usepackage[a4paper, left=1.0cm, right=1.0cm, top=2.5cm, bottom=2.5cm]{geometry}

\usepackage{amsmath}
\usepackage{amssymb}
\usepackage{xspace}
\usepackage{soul}
\usepackage{mathtools}
\usepackage{makecell}
\usepackage{cite}
\usepackage{nicefrac}
\usepackage{float}
\usepackage{wrapfig}
\usepackage{amsthm}

\newtheorem{definition}{Definition}
\newtheorem{remark}{Remark}
\newtheorem{proposition}{Proposition}
\newtheorem{lemma}{Lemma}
\newtheorem{corollary}{Corollary}
\newtheorem{theorem}{Theorem}

\newcommand{\ee}{\varepsilon}
\newcommand{\N}{\mathbb{N}}

\newcommand{\R}{\mathbb{R}}
\newcommand{\UU}{\mathbb{U}}
\newcommand{\DD}{\mathbb{D}}
\newcommand{\II}{\mathbb{I}}

\newcommand{\LL}{\mathcal{L}}
\newcommand{\U}{\mathcal{U}}
\newcommand{\G}{\mathcal{G}}
\newcommand{\D}{\mathcal{D}}
\newcommand{\RA}{\mathcal{R}}
\newcommand{\RAT}{\mathcal{R}_{T}}
\newcommand{\RAC}{\mathcal{R}^{\mathsf{C}}}
\newcommand{\RATC}{\mathcal{R}^{\mathsf{C}}_{T}}
\newcommand{\RAB}{\partial \RA}
\newcommand{\RABT}{\partial \RAT}
\newcommand{\RAM}{[\partial \RA]_-}
\newcommand{\RAO}{[\partial \RA]_0}

\newcommand{\cl}[1]{\mathsf{cl}(#1)}
\newcommand{\Int}[1]{\mathsf{int}(#1)}

\newcommand{\ie}{i.e.~}     
\newcommand{\eg}{e.g.~}

\newcommand{\keywords}[1]{\noindent\textbf{Keywords:} #1}

\begin{document}

\title{On the Boundary of the Robust Admissible Set in State and Input Constrained Nonlinear Systems}

\author[1]{Franz Ru{\ss}wurm}

\author[2]{Jean Lévine}

\author[1]{Stefan Streif\thanks{Corresponding author: stefan.streif@etit.tu-chemnitz.de}}

\affil[1]{Automatic Control and System Dynamics, University of Technology Chemnitz, Germany}

\affil[2]{Unité Maths et Systèmes, MINES Paris - PSL University, France}

\date{}
\maketitle

\abstract{In this paper, we consider nonlinear control systems subject to bounded disturbances and to both state and input constraints. We introduce the definition of \emph{robust admissible set} - the set of all initial states from which the state and input constraints can be satisfied for all times against all admissible disturbances. We focus on its boundary that can be decomposed into the \emph{usable part} on the state constraint boundary and the \emph{barrier}, interior to the state constraints. We show that, at the intersection of these two components, the boundary of the robust admissible set must be tangent to the state constraint set and separate the interior of the robust admissible set and its complement, a property that we call the \emph{ultimate locally separating hyperplane condition}. Moreover, we prove that the barrier must satisfy a \emph{saddle-point principle} on a Hamiltonian, based on Pontryagin's maximum principle, whose final condition is precisely the ultimate locally separating condition, thus providing a set of differential equations made of the system and its adjoint for a direct construction of the barrier. Lastly, we illustrate our results by calculating the robust admissible set for an adaptive cruise control example.}


\renewcommand\thefootnote{}

\renewcommand\thefootnote{\arabic{footnote}}

\section{Introduction}\label{Sec:Introduction}
In this paper, we consider nonlinear control systems, which are affected by disturbances and submitted to state and input constraints.
For these systems, we aim at characterizing the \emph{robust admissible set}, \ie the set of all initial states for which a control input exists, such that the associated integral curve satisfies the state and input constraints for all times whatever the disturbances are.

Our approach extends the one presented in \cite{Levine_2013} for ordinary (single-valued) differential systems with disturbances, a version of the \emph{viability kernel} theory introduced by \cite{Aubin_1991} in the context of multi-valued differential systems.
The so-called \emph{admissible set} of \cite{Levine_2013} becomes here the \emph{robust admissible set} and the part of its boundary running in the interior of the constraint set is still called the \emph{barrier}. 

The concept of barrier was introduced by Isaacs \cite{Isaacs_1965} in the context of differential games. It plays here a central role not only because it separates the admissible states from the non admissible ones, but also because it constitutes a \emph{semi-permeable surface}, namely a surface that can be one-way crossed, without possible return.

Successful applications have been made to systems without disturbances in several fields such as  \emph{food production} systems \cite{Aschenbruck_2020}, \emph{power grids} \cite{Aschenbruck_2021} and \emph{maintaining infection caps during epidemics} \cite{Esterhuizen_2021}, as well as some problems of Model Predictive Control \cite{Esterhuizen_et-aL_2021}, to cite just a few ones.

Let us first illustrate the topic of the present paper by the following example of adaptive cruise control.
It was originally introduced in a different context in \cite{Ames_2014} and then extended to include disturbances in \cite{Xu_2015}.
It is a simple model of two cars driving in a convoy, constrained to keep a safe distance between each other.

\subsection{Introductory Example: Adaptive Cruise Control}\label{intro-ex:subsec}
We consider an adaptive cruise control scenario involving two vehicles driving in a convoy, the leading one, whose speed is denoted by $x_{1}$ and the follower, whose speed is denoted by $x_{2}$. The distance between them is denoted by $x_{3}$.

The triple $(x_{1}, x_{2}, x_{3})$, the system's state, is supposed to satisfy the dynamical system
\begin{align}\label{eq:ACC}
\begin{split}
\dot{x}_1(t) &= a + d_1(t), \\
\dot{x}_2(t) &= -\left(a_0 + a_1 x_2(t) + a_2 x_2^2(t)\right) + g d_2(t) + g u(t), \\
\dot{x}_3(t) &= x_1(t) - x_2(t),
\end{split}
\end{align}
where 
\begin{itemize}
\item $d_{1}$ is an unknown disturbance of the leader's speed $a$, known to the follower, and $d_1(t)$ is subject to the constraint 
$$
d_1(t) \in \DD_1 \triangleq [\underline{d}_1,\overline{d}_1] \subset \R.
$$
\item  the force applied to the follower is equal to the sum of the wheel force, the acceleration $u$ (the control) multiplied by its mass $m$ (in kg), and the aerodynamic drag and rolling resistance $-(a_0 + a_1 x_2(t) + a_2 x_2^2(t)) + g d_2(t)$ (in N), $d_{2}(t)$ being a second unknown disturbance subject to  
$$
d_2(t) \in \DD_2 \triangleq [\underline{d}_2,\overline{d}_2] \subset \R,
$$
and $g$ the gravitational force $g = 9.81 \textnormal{ms}^{-2}$.
\end{itemize}
The control $u$ is assumed to be limited by its maximal acceleration $\overline{u}$ and deceleration  $\underline{u}$:
\begin{align*}
    u\in \UU \triangleq [\underline{u},\overline{u}] \subset \R.
\end{align*} 
We impose the following state constraints for the sake of safety
\begin{align}\label{eq:ACC_constraints}
\begin{split}
g_1(x) &\triangleq \tau x_2 - x_3 \leq 0, \\
g_2(x) &\triangleq x_3 - D_{\textnormal{max}} \leq 0,
\end{split}
\end{align}
where $\tau$ is the minimum time headway, and $D_{\textnormal{max}}$ the maximum allowed distance between the two vehicles.

Determining the set of initial states for which there exists a control of the follower $t\mapsto u(t) \in \UU$ for all $t$,  such that, despite the presence of the unknown disturbances $t\mapsto d_{1}(t)\in \DD_1$ and $t\mapsto d_{2}(t)\in \DD_2$, the integral curves of System~\eqref{eq:ACC} satisfy the safety constraints $g_1$ and $g_{2}$ for all future times, constitutes the archetypal robust admissible problem that we address in this paper and solve in Section~\ref{Sec:Example}.

\subsection{Some Historical Overview of the Problem}
In the case without disturbance the admissible set has been thoroughly studied in \cite{Levine_2013} and then in \cite{Esterhuizen_2017,WE-JL_automatica_2016} for mixed constraints, \ie constraints depending of both the state and the controls. 
Following the approach of \cite{Levine_2013}, an extension to uncontrolled systems subject to disturbances has been done in \cite{Esterhuizen_2020}, where the notion of \emph{maximal robust positively invariant set} is introduced. 
The importance of these approaches relies on the fact that the admissible set is easily described by its boundary. The latter is partly constituted by the so-called \emph{barrier}, made of specific system trajectories satisfying a minimum principle and enjoying the \emph{semi-permeability property}, introduced by Isaacs  \cite{Isaacs_1965} in the context of \emph{differential games}.

Other approaches to study robust admissible trajectories satisfying the constraints for continuous time nonlinear systems include  \emph{Hamilton-Jacobi backwards reachability} \cite{Mitchell_2005,Cruck_2008,Chen_2018,Herbert_2021,Maidens_2013}.
This problem, interpreted as a \emph{differential game} \cite{Isaacs_1965} with two opposing inputs, is usually solved by considering \emph{viscosity solutions} of the Hamilton-Jacobi-Bellman partial differential equation.

The admissible set also plays an important role in model-predictive control (MPC) \cite{Rawlings_2017,Esterhuizen_et-aL_2021,Santos_2021}, with or without disturbances.
It has been remarked in \cite{Kim_2021} that it can be used to enlarge the terminal region and improve stability performance.

The admissible set being the largest invariant set in some sense, let us briefly mention approaches that have been developed in the context of \emph{forward reachable sets} \cite{Sontag_1998}, \emph{invariant sets} \cite{Blanchini_2015}, geometric deformations of polytopes \cite{Ameen_2023}, and other related concepts as the notion of \emph{capture basin}, also called \emph{attainable region}, the set of all initial states for which a target can be reached in finite time \cite{Feinberg_2002}.
Similar works focused on \emph{inner and outer approximations} based on \emph{interval analysis}, \cite{Monnet_2016,Martin_2018}.

Finally, \emph{Control barrier functions} (CBFs) \cite{Choi_2021,Ames_2019}, have received significant attention in recent years for safety critical control, and a solution of the above introductory example has been proposed in this context (see \cite{Ames_2014,Xu_2015}), only with the constraint $g_1$ in \eqref{eq:ACC_constraints}. Though both, CBFs and barriers of admissible sets wear similar names, and aim at ensuring constraint satisfaction, they differ in how safety is enforced: CBFs are \emph{Control-Lyapunov functions} that render a given set forward invariant by feedback  \cite{Ames_2017,Xu_2018,Zhang_2023,Jankovic_2018,Clark_2021}, whereas the barrier of the robust admissible set exactly characterizes the boundary of the maximal set of initial conditions from which constraint satisfaction can be guaranteed for all future times, without need to design any feedback control law.

\subsection{Content and Contributions of the Paper}
In this paper, we extend the characterization of the admissible set's boundary from \cite{Levine_2013} to nonlinear control systems subject to bounded disturbances, using a comparable mixed geometric and topological approach, where we derive a \textit{saddle point} version of the Pontryagin maximum principle, that holds on a part of the boundary of the robust admissible set.
More precisely, we show that this condition holds at points where the boundary of the robust admissible set intersects tangentially with the active state constraint and locally separates the state space between the robust admissible set and its complement. From such points it extends to extremal trajectories evolving towards the interior of the constrained state space, obtained via the solution of an adjoint differential equation and the control and disturbance must satisfy a saddle point condition of a Hamiltonian.

Most importantly with our approach, solving a partial differential equation such as the Hamilton-Jacobi-Bellman equation is not needed and, in contrast to iterative approximation methods, such as \eg interval analysis, we obtain the boundary of the robust admissible set directly by integration of an ordinary differential equation, without need of inner or outer approximation.

Let us insist on the fact that, while Hamilton–Jacobi–Bellman and reachability methods typically rely on the computation of a value function over the full state space and of the corresponding optimal trajectories in feedback form, the barrier approach focuses directly on the $(n-1)$-dimensional boundary of the admissible set, thereby reducing the dimensionality of the problem. Note in addition that this boundary is obtained by solving the local optimization problem comprised in the minimum principle \emph{without state constraints}, unlike other approaches based on the afore-mentioned Hamilton-Jacobi-Bellman equation or interval analysis. Hence, our approach opens some modest but interesting new perspectives to tackle higher dimensional problems \cite{Esterhuizen_2015}.

Moreover, if the problem depends on parameters, the boundary may be obtained as a function of these parameters, without requiring to re-compute the whole admissible set for each parameter value, as it is sometimes needed with numerical methods. 

\bigskip

Our main contributions are as follows:
\begin{itemize}
    \item We characterize the geometry of the boundary of the robust admissible set. Our first contribution cpncerns the extension of the theory developed in \cite{Levine_2013} for nonlinear control systems and in \cite{Esterhuizen_2020} for disturbed but uncontrolled nonlinear systems, to the setting of disturbed and controlled nonlinear systems subject to state constraints. In particular, we prove that the boundary must intersect the active state constraints tangentially and that the corresponding tangent hyperplane locally separates the state space between the robust admissible set and its complement.
    \item We derive a necessary condition in the form of a saddle point condition of the Hamiltonian on the barrier characterizing optimal control and worst-case disturbance inputs, which constitutes our second contribution. Let us stress that in spite of the min-max natural structure of the problem (worst-case design), the obtained conditions take in fact the form of a saddle point condition on the Hamiltonian (not only a min-max) at every point of the boundary.
    \item We present an application to the adaptive cruise control problem of the introduction that has been studied in \cite{Ames_2014,Xu_2015} via the approach of control barrier functions. Our contribution is the complete and explicit description of the largest robust admissible set $\RA$ through its boundary $\RAB$, parameterized by the leading car's speed. It should be remarked that our approach is not limited to a single constraint as in \cite{Ames_2014,Xu_2015}.    
The control barrier function framework appears to be complementary to ours in the sense  that it is able to design a feedback that ensures  the asymptotic convergence of the follower to a subset $\cal{C}$ of $\RA$, given some disturbances.
\end{itemize}

\bigskip

The outline of this paper is as follows.

In Section \ref{Sec:The_Robust_Admissible_Set}, the control system, the constraints and the necessary assumptions are formally introduced.
Moreover, a definition of the robust admissible set in finite and infinite time is given, together with a proposition establishing its closedness, hence  including its boundary.
Section \ref{Sec:The_Boundary_of_the_Robust_Admissible_Set} develops analytical tools to characterize these boundary parts and in particular the so-called \emph{barrier}, as well as the associated barrier trajectories.
Section \ref{Sec:Ultimate_Sep} presents our main result characterizing the barrier trajectories in terms of adjoint dynamics, Hamiltonian and  a saddle point type necessary condition.
In Section \ref{Sec:Example}, we apply the developed method to the adaptive cruise control system introduced in Section~\ref{Sec:Introduction}.
Finally, Section \ref{Sec:Conclusion} concludes the paper.

\subsection{Notation}
Throughout the paper, we use the following notations:
\begin{itemize} 
\item The symbol $\triangleq$ is used to indicate an identity by definition.
\item For a scalar $x \in \R$, $\left\vert x \right\vert$ denotes its absolute value, while $\Vert x \Vert$ denotes the Euclidean norm of a vector $x \in \R^n$ and $x^\top$ its transpose.
\item Given $p\in \N$ and $g\triangleq (g_{1}, \ldots, g_{p})^{T}$ a vector of real valued functions $g_{i}$, $i=1, \ldots,p$, the notation $g(x)\preceq 0$ (resp. $g(x)\prec 0$) means that $g_{i}(x)\leq 0$ (resp. $g_{i}(x) < 0$) for all $i= 1, \ldots, p$.
\item For sets $X$ and $Y$, $\mathcal{C}^k(X,Y)$ denotes the space of $k$-times continuously differentiable functions mapping $X$ to $Y$.
\item For a set $S$, $S^C$, $\Int{S}$, $\cl{S}$ and $\partial S$ denote, respectively, its complement, interior, closure and boundary.
\end{itemize}

\section{The Robust Admissible Set}\label{Sec:The_Robust_Admissible_Set}

Consider a nonlinear control system with two kinds of inputs denoted by $u$, the controls, and $d$, the disturbances, \ie  the non controlled input variables,
\begin{align}
\begin{split}\label{eq:control_system}
\dot{x}(t) &= f(x(t),u(t),d(t)) \\
x(t_0) &= x_0 \\
u &\in \U, \; \; d \in \D
\end{split} 
\end{align}
with state $x(t) \in \R^n$, the set $\U$ (resp.$\D$)  of open-loop control laws (resp. disturbances) consist of Lebesgue-measurable functions $u : [t_0,\infty) \rightarrow  \UU$ with $\UU \subset \R^m$ nonempty and compact (resp. $d:[t_0,\infty) \rightarrow \DD$  with $\DD \subset \R^w$ nonempty and compact).

We denote by $x^{(x_0,u,d)}$ the absolutely continuous maximal integral curve that satisfies \eqref{eq:control_system}, initiated at $x_0 \in \R^n$ and generated by a control law $u \in \U$ and a disturbance $d \in \D$. 

In addition, given an arbitrary natural number $p$, we consider state constraints of the form
\begin{align} 
g_i(x(t)) \leq 0 \; \; \; \; \forall \, t \in [t_0, \infty), \; \; \forall \; i \in \{1,2, \ldots, p \}.\label{eq:state_constraints}
\end{align}
We denote by $g(x)$ the vector $(g_1(x),g_2(x),\ldots,g_p(x))^\top$ and define the constraint set $\G  \triangleq \{x \in \R^n \, \vert \, g(x) \preceq 0 \}$, which is assumed to have a nonempty interior $\G_- \triangleq \{x \in \R^n \, \vert \, g(x) \prec 0 \}$ and and we denote its boundary by $\G_0 \triangleq \{x \in \R^n \, \vert \, \exists \, i \in \{1,2,\ldots,p \}: g_i(x) = 0 \}$, assumed to be non empty. We indeed have $\G = \G_0 \cup \G_-$.

Following \cite{Levine_2013}, we make the assumptions:
\begin{description}
	\item[Assumptions]
\end{description}
\begin{enumerate}
    \item[(H1)] $f$ is an at least $\mathcal{C}^2$ vector field of $\R^n$, whose dependence with respect to $u$ and $d$ is at least $\mathcal{C}^2$ in an open subset of $\R^m \times \R^w$ containing $\UU \times \DD$.
    \item[(H2)] There exists a constant $0<C<\infty$ such that
    \begin{align*}
        \sup_{u \in \UU,~ d\in  \DD} \left\vert x^\top f(x,u,d) \right\vert \leq C (1 + \Vert x \Vert^2) \; \; \forall \, x\in\R^n.
    \end{align*}
   
    \item[(H3)] The set $f(x,\UU,\DD) \triangleq \{ f(x,u,d)  \mid u \in \UU, \, d\in \DD \}$, called the vectogram, is convex for all $x\in\R^n$.
    \item[(H4)] For all $i \in \{1,2,\ldots,p \}$, $g_i: \R^n \rightarrow \R$ is an at least $\mathcal{C}^2$ function and the set $\left\{x \in \R^n \, \big\vert \, \max_{i \in \{1,2,\ldots,p \}} g_i(x) = 0 \right\}$ defines a piecewise  $\mathcal{C}^1$ manifold of dimension $n-1$.
\end{enumerate}

\begin{remark}\label{Rem:Assumptions}
Assumptions (H1) and (H2) ensure the existence, uniqueness and boundedness on $[t_{0}, \infty)$ of the integral curves of \eqref{eq:control_system} for all $u \in \U$ and $d \in \D$. Assumption (H2), which means that $f$ grows at most linearly with respect to $x$, more precisely that $\sup_{(u,d)\in \UU\times\DD} ||f(x, u,d)|| \leq C(1+||x||)$ for all $x\in \R^n$, as discussed in \cite{Levine_2013}, is a sufficient condition (but not necessary) to guarantee the boundedness and compactness of the set of integral curves used in the analysis. 
\\
In many applications, similar properties of the set of integral curves can be ensured under weaker but less general assumptions, for instance by restricting the constraint set to a compact region of the state space, as in Example~\ref{intro-ex:subsec}. The interested reader is invited to consult the discussion on this subject \eg in \cite{Filippov_siam} or in the book \cite[Chapter 10, section 10.3]{Agrachev_2004}.

Assumption (H3) holds, in particular in the class of systems that are affine in the control and disturbance inputs, \ie 
$f(x, u, d) = f_{0}(x) + u~f_{1}(x) + d~f_{2}(x)$, when $\UU$ and $\DD$ are convex.

Finally, assumption (H4) is a natural regularity assumption meanning that the set
\begin{align*}
\left\{x \in \R^n \, \big\vert \, \max_{i \in \{1,2,\ldots,p \}} g_i(x) = 0 \right\}
\end{align*}
admits almost everywhere a tangent $n-1$ dimensional plane and a 1-dimensional non degenerate normal.
\end{remark}

The following compactness result is a fundamental tool used in most of the results of this paper, and in particular in the fact that the boundary of the robust admissible set, defined on the inifinite interval of time $[0, \infty)$, is well defined (see Proposition~\ref{Prop:RAT_closed}), as well as its consequences on  Propositions~\ref{Prop:Characteriztion_RAB}, \ref{Prop:RAB_max_equal_0},~\ref{Prop:Characterizing_Trajectories_in_RAM},~\ref{Prop:Ultimate_Separation}). 
\begin{proposition}\label{PropA2}
Let us assume that (H1)--(H3) are valid.
 Given a compact subset $X_0 \subset \R^n$, the set of all integral curves satisfying \eqref{eq:control_system} and \eqref{eq:state_constraints}, with initial values $x_0 \in X_0$, is compact with respect to the topology of uniform convergence on $\mathcal{C}^0([0,T],\R^n)$ for all $T>0$.
\end{proposition}
\begin{proof}
The proof follows the exactly same lines as in \cite[Corollary A.1 and Lemma A.2]{Levine_2013} with the input $u \in \U$ replaced by the pair $v=(u,d) \in \U\times \D$.
\end{proof}

We next define the robust admissible set of system \eqref{eq:control_system} subject to the state constraints \eqref{eq:state_constraints}.
\begin{definition}\label{Def:robust_admissible_set}
The \emph{robust admissible set}, denoted by $\RA$, is the set of all initial states $\bar{x}$ in the constraint  set $\G$ for which there exists a control $u\in \U$ such that, for all disturbance $d\in  \D$,  the corresponding integral curve $x^{(\bar{x},u,d)}$ of system \eqref{eq:control_system} starting at $\bar{x}$,  satisfies the state constraints \eqref{eq:state_constraints} for all future times. In mathematical terms it reads: 
\begin{equation}\label{RA:def}
\RA \triangleq \{ \bar{x} \in \G \mid \exists  u \in \U : x^{(\bar{x},u,d)}(t) \in \G ~ \forall d \in \D, \, \forall t \geq 0 \} .
\end{equation}
\end{definition}

Its complement  relative to $\G$ is given by
\begin{align*}
    \RAC \triangleq \G \setminus \RA = \{ \bar{x} \in \G \mid \forall u \in \U, \exists d \in \D, \, \exists \bar{t} < \infty: ~x^{(\bar{x},u,d)}(\bar{t}) \not\in \G \} .
\end{align*}

\begin{remark}
The determination of the robust admissible set \eqref{RA:def}, might also be viewed as a \emph{min-max problem} in the framework of qualitative differential games with perfect information \cite{Isaacs_1965}, where $u$ is played by player 1, the minimizer, and $d$ by player 2, the maximizer. The set $\RA$ could  be relabeled $\RA^-$ to distinguish it from its \emph{max-min} counterpart $\RA^+$: 
\begin{align*}
\RA^+ \triangleq \{ \bar{x} \in \G \mid \forall d \in \D, \, \exists u \in \tilde{\U} : ~x^{(\bar{x},u,d)}(t) \in \G \; \; \forall  t \geq 0 \},
\end{align*}
where $\tilde{\U}$ is the set of non-anticipative functions $u: (t,d) \in[0, \infty)\times \DD \mapsto u(t,d)\in \UU $, regular enough to ensure existence and uniqueness of $x^{(\bar{x},u,d)}$, \ie $t\mapsto u(t,d(t))$ measurable for all $d\in \D$. It means that  the disturbance $d$ must be measured and communicated to player 1 at all times, so that player 1 can instantaneously react to $d(t)$ by choosing a suitable $u(t))$. We indeed have $\U \subset \tilde{\U}$ and thus the minimum over $\U$ may be larger than or equal to the minimum over $\tilde{\U}$.
However, since $t\mapsto u\circ d(t) \triangleq u(t,d(t))$ is measurable for all $d\in \D$, we also deduce that $u\circ d \in \U$, hence proving that the minimum over $\U$ and $\tilde{\U}$ are equal. 
In this paper, we focus on the study of $\RA$ given by \eqref{RA:def}, which means that the minimizer does not acquire information during the game on the maximizer's input and symmetrically, the maximizer does not acquire information during the game on the minimizer's input.\footnote{We do not discuss here information structures where the players have access to exact or noisy observations, functions of the state. Such discussion is far beyond the scope of this paper.}
\end{remark}

We also consider the {finite horizon robust admissible set}
\begin{align*}
\RAT \triangleq \{ \bar{x} \in \G \, \vert \, \exists \, u \in \U : \; x^{(\bar{x},u,d)}(t) \in \G \; \forall \, d \in \D, \, \forall \, t \in [0, T] \}
\end{align*}
as well as its complement relative to $\G$
\begin{align*}
    \RATC \triangleq \G \setminus \RAT = \{ \bar{x} \in \G \, \vert \, \forall u \in \U, \, \exists \, d \in \D, \, \exists \, \bar{t} \in [0,T]: \; x^{(\bar{x},u,d)}(\bar{t}) \not\in \G \}.
\end{align*}

\begin{proposition}\label{Prop:RAT_closed}
Under assumptions (H1)-(H4), $\RAT$ is closed for all finite $T>0$, as well as $\RA = \underset{T \geq 0}{ \bigcap } \RAT$, and thus contain their respective boundary $\RABT$ and $\partial \RA$.
\end{proposition}
\begin{proof}
The proof is adapted from the one of \cite[Proposition 4.1]{Levine_2013} as follows.
Consider a sequence of initial states $\{ x_{k} \}_{k\in \N}$ in $\RAT$ converging to $\bar{x}$ as $k$ tends to infinity. By definition of $\RAT$, for every $k\in \N$, there exists $u_{k}\in \U$ such that the corresponding integral curve $x^{(x_{k},u_{k},d)}$ satisfies $g(x^{(x_{k},u_{k},d)}(t)) \preceq 0$ for all $d\in \D$ and $t\in [0,T]$. According to Proposition~\ref{PropA2}, there exists a uniformly converging subsequence, still denoted by $x^{(x_{k},u_{k},d)}$, to the absolutely continuous integral curve $x^{(\bar{x},\bar{u},d)}$ for some $\bar{u}\in \U$ and for all $d\in \D$. By the continuity of $g$, we immediately get that
$g(x^{(\bar{x},\bar{u},d)}(t)) \preceq 0$ for all $d\in \D$ and all $t\in [0,T]$, hence $\bar{x}\in \RAT$, and the closedness of $\RABT$ is proven.

The closedness of $\RA$ follows from the fact that, for all $0\leq T_{1} \leq T_{2} < \infty$, we have $\RA = \RA_{\infty} \subset \RA_{T_{2}} \subset \RA_{T_{1}} \subset \RA_{0} = G$, thus
$\RA = \underset{T \geq 0}{ \bigcap } \RAT$, and the intersection of a family of closed sets is closed.
\end{proof}

\section{The Boundary of the Robust Admissible Set}\label{Sec:The_Boundary_of_the_Robust_Admissible_Set}

We now proceed in characterizing $\RA$, $\RAC$ and the boundary  $\RAB$. 
To this aim, we introduce a qualification assumption to guarantee that no gap exists between the corresponding criteria:
\begin{itemize}
\item[(H5)] There exists an $\bar{x} \in \RAB$ and a $u\in \U$ such that, for every $d \in \D$, $\sup_{t \in [0,\infty)} \, \max_{i = 1,2, \ldots,p } \, g_i(x^{(\bar{x},u,d)}(t)) = 0$.
\end{itemize}
\begin{remark}This assumption, that was forgotten in \cite{Levine_2013}, may be interpreted as a robust-output-null-controllability condition.

Note that the function $(\bar{x},t) \mapsto \max_{i = 1,2, \ldots,p } g_i(x^{(\bar{x},u,d)}(t))$ is continuous for every $u$ and $d$ but that, taking the supremum with respect to $t$, downgrades its continuity to upper semi-continuity with respect to $\bar{x}$ due to a possible jump over the 0 value when $\bar{x}$ tends to a boundary point of the constraint set while remaining in its interior, hence the assumption (H5) to exclude this case.
\end{remark}

\begin{proposition}\label{Prop:Characteriztion_RAB}
    Assume that assumptions (H1)-(H4) hold.
    Then:
    \begin{enumerate}
        \item[(i)] $\bar{x} \in \RA$ is equivalent to
        \begin{align}\label{eq:RA_leq_0}
        \min_{u \in \U} \, \sup_{d \in \D} \, \sup_{t \in [0,\infty)} \, \max_{i = 1,2, \ldots,p } \, g_i(x^{(\bar{x},u,d)}(t)) \leq 0 .
        \end{align}
        \item[(ii)] $\bar{x} \in \RAC$ is equivalent to
        \begin{align}\label{eq:RAC_greater_0}
        \min_{u \in \U} \, \sup_{d \in \D} \, \sup_{t \in [0,\infty)} \, \max_{i = 1,2, \ldots,p } \, g_i(x^{(\bar{x},u,d)}(t)) > 0 .
        \end{align}
        \item[(iii)] If (H5) holds, then $\bar{x} \in \RAB$ is equivalent to
        \begin{align}\label{eq:RAB_equal_0}
        \min_{u \in \U} \, \sup_{d \in \D} \, \sup_{t \in [0,\infty)} \, \max_{i = 1,2, \ldots,p } \, g_i(x^{(\bar{x},u,d)}(t)) = 0 .
        \end{align}
    \end{enumerate}
\end{proposition}
\begin{proof}
    We start by proving (i).
    Assume $\bar{x} \in \RA$.
    By definition of $\RA$, there exists a control law $u \in \mathcal{U}$ such that $x^{(\bar{x},u,d)}(t) \in \RA$ for all $t>0$ and all disturbances $d \in \mathcal{D}$.
    We get
    \begin{align*}
        \sup_{d \in \mathcal{D}} \, \sup_{t \in [0,\infty)} \, \max_{i = 1,2, \ldots,p } \, g_i(x^{(\bar{x},u,d)}(t)) \leq 0
    \end{align*}
    and, by the definition of the infimum, we immediately get
    \begin{align}\label{eq:proof_inf_u_sup_d_leq_0}
        \inf_{u \in \mathcal{U}} \, \sup_{d \in \mathcal{D}} \, \sup_{t \in [0,\infty)} \, \max_{i = 1,2, \ldots,p } \, g_i(x^{(\bar{x},u,d)}(t)) \leq 0.
    \end{align}
    We proceed by proving that the infimum is achieved by some $\bar{u} \in \mathcal{U}$.
    To this end, we consider a sequence $u_k \in \mathcal{U}$ such that
    \begin{equation}\label{eq:proof_limit_g(x_k)}
		\lim_{k \rightarrow \infty} \, \sup_{d \in \mathcal{D}} \, \sup_{t \in [0, \infty)} \, \max_{i = 1,2,\ldots,p} g_i(x^{(\bar{x},u_k,d)}(t))
		=   \inf_{u \in \mathcal{U}} \, \sup_{d \in \mathcal{D}} \, \sup_{t \in [0,\infty)} \, \max_{i = 1,2, \ldots,p } \, g_i(x^{(\bar{x},u,d)}(t)).
    \end{equation}
    As shown in the proof of Proposition \ref{Prop:RAT_closed} and by taking $x_k = \bar{x}$ for all $k \in \mathbb{N}$, there exists a  converging subsequence $x^{(\bar{x},u_{k_l},d)}$ uniformly with respect to $d$,  on every compact interval $[0,T]$ with $T>0$, whose limit is $x^{(\bar{x},\bar{u},d)}$  for all $d \in \mathcal{D}$ and for some $\bar{u} \in \mathcal{U}$, independent of $d$.
    From the continuity of $g$ we also obtain that  the sequence $g_i(x^{(\bar{x},u_{k_l},d)})$ uniformly converges to $g_i(x^{(\bar{x},\bar{u},d)})$ as $l \rightarrow \infty$ for all $d \in \mathcal{D}$ on every compact interval $[0,T]$ with $T>0$ and for all $i = 1,2, \ldots,p$.
    Thus, for every $T>0$ and every $\varepsilon > 0$, there exists an $l_0(T,\varepsilon) \in \mathbb{N}$ such that for all $l \geq l_0(T,\varepsilon)$ and all $i = 1,2,\ldots,p$,
    \begin{align*}
        \left\vert g_i(x^{(\bar{x},\bar{u},d)}(t)) - g_i(x^{(\bar{x},u_{k_l},d)}(t)) \right\vert \leq \varepsilon
    \end{align*}
   holds for all $t \in [0,T]$ and all $d \in \mathcal{D}$.
    By taking the supremum w.r.t. $d \in \mathcal{D}$ and $t \in [0,\infty)$, as well as the maximum over all $i = 1,2,\ldots,p$, we obtain
    \begin{align*}
        g_i(x^{(\bar{x},\bar{u},d)}(t)) \leq \sup_{d \in \mathcal{D}} \, \sup_{t \in [0,T)} \, \max_{i = 1,2,\ldots,p} \, g_i(x^{(\bar{x},u_{k_l},d)}(t)) + \varepsilon
    \end{align*}
    for all $t \in [0,T]$ and $d \in \mathcal{D}$.
    On the other hand, according to the definition of the limit \eqref{eq:proof_limit_g(x_k)} and taking into account the uniform convergence of the above-mentioned subsequence, we obtain that for every $\varepsilon > 0$ there exists a $l_1(\varepsilon) \in \mathbb{N}$ such that, for all $l \geq l_1(\varepsilon)$,
   $$
    		\sup_{d \in \mathcal{D}} \, \sup_{t \in [0,\infty)} \, \max_{i = 1,2,\ldots, p} \, g_i(x^{(\bar{x},u_{k_l},d)}(t)) \leq \inf_{u \in \mathcal{U}} \, \sup_{d \in \mathcal{D}} \, \sup_{t \in [0,\infty)} \, \max_{i = 1,2,\ldots, p} \, g_i(x^{(\bar{x},u,d)}(t)) + \varepsilon.
  $$
    We conclude that, for all $l \geq \max \{l_0(T,\varepsilon), l_1(\varepsilon) \}$,
    \begin{align*}
        g_i(x^{(\bar{x},\bar{u},d)}(t)) \leq \inf_{u \in \mathcal{U}} \, \sup_{d \in \mathcal{D}} \, \sup_{t \in [0,\infty)} \, \max_{i = 1,2,\ldots, p} \, g_i(x^{(\bar{x},u,d)}(t)) + 2 \varepsilon.
    \end{align*}
    for all $t \in [0,T]$, all $d \in \mathcal{D}$ and for all $i = 1,2,\ldots,p$.
    This inequality does not depend on $k_l$ and holds for all $T>0$.
    Moreover, the right-hand side is independent of the index $i$, the time $t$, the time horizon $T$ and the disturbance $d$.
    Thus, we can take the supremum w.r.t. $t \in [0,\infty)$ and $d \in \mathcal{D}$ as well as the maximum over $i = 1,2,\ldots,p$ of the left-hand side.
   Finally, according to the definition of the infimum w.r.t. $u \in \mathcal{U}$, we obtain
\begin{align*}
\sup_{d \in \mathcal{D}} \, \sup_{t \in [0,\infty)} \, \max_{i = 1,2,\ldots, p} \, g_i(x^{(\bar{x},\bar{u},d)}(t))  &\leq \inf_{u \in \mathcal{U}} \, \sup_{d \in \mathcal{D}} \, \sup_{t \in [0,\infty)} \, \max_{i = 1,2,\ldots, p} g_i(x^{(\bar{x},u,d)}(t)) + 2 \varepsilon \\
&\leq \sup_{d \in  \mathcal{D}} \, \sup_{t \in [0,\infty)} \, \max_{i = 1,2,\ldots, p} \, g_i(x^{(\bar{x},\bar{u},d)}(t)) + 2 \varepsilon
\end{align*}
    for all $\varepsilon > 0$, hence the equality and the minimum over $u \in \mathcal{U}$ is achieved at some $\bar{u}\in \mathcal{U}$, which, together with \eqref{eq:proof_inf_u_sup_d_leq_0}, implies \eqref{eq:RA_leq_0}.

    Conversely, assume that \eqref{eq:RA_leq_0} holds.
    Then, there exists a control input $u \in \mathcal{U}$ such that, for all disturbances $d \in \mathcal{D}$, we have
    \begin{align*}
        \sup_{t \in [0,\infty)} \, \max_{i = 1,2,\ldots, p} \, g_i(x^{(\bar{x},u,d)}(t)) \leq 0.
    \end{align*}
    Thus, $g_i(x^{(\bar{x},u,d)}(t)) \leq 0$ for all $t \geq 0$, all $d \in \mathcal{D}$ and all $i = 1,2,\ldots,p$, which proves that $\bar{x} \in \RA$.

    In order to prove (ii), we assume that $\bar{x} \in \RAC$.
    From the definition of $\RAC$, there exists a disturbance $d \in \mathcal{D}$ such that, for all $u \in \mathcal{U}$,
    \begin{align*}
        \sup_{t \in [0,\infty)} \, \max_{i = 1,2,\ldots,p} \, g_i(x^{(\bar{x},u,d)}(t)) > 0.
    \end{align*}
    Thus
    \begin{align*}
        \inf_{u \in \mathcal{U}} \, \sup_{d \in \mathcal{D}} \, \sup_{t \in [0,\infty)} \, \max_{i = 1,2,\ldots,p} \, g_i(x^{(\bar{x},u,d)}(t)) \geq 0.
    \end{align*}
    We now resume the argument of uniformly converging sequence as in the proof of (i). Hence, the minimum over $u \in \mathcal{U}$ is achieved by some  $\bar{u} \in \mathcal{U}$ and we have
    \begin{align*}
        \min_{u \in \mathcal{U}} \, \sup_{d \in \mathcal{D}} \, \sup_{t \in [0,\infty)} \, \max_{i = 1,2,\ldots,p} \, g_i(x^{(\bar{x},u,d)}(t)) \geq 0.
    \end{align*}
    Note that this inequality has to be strict.
    Otherwise, the equality to 0 would imply that $\bar{x} \in \RA$, thus contradicting that $\bar{x} \in \RAC$ by assumption.
    
    Conversely, assume that \eqref{eq:RAC_greater_0} holds.
    Thus for all $u \in \mathcal{U}$ there exists a disturbance $d \in \mathcal{D}$ such that
    \begin{align*}
        \sup_{t \in [0,\infty)} \, \max_{i = 1,2,\ldots,p} \, g_i(x^{(\bar{x},u,d)}(t)) > 0.
    \end{align*}
    The continuity of the mapping $t \longmapsto \max_{i = 1,2,\ldots,p} \, g_i(x^{(\bar{x},u,d)}(t))$ thus implies that the inverse image 
    \begin{align*}
        \left\lbrace t \in [0,\infty) \; \Big\vert \; \max_{i = 1,2,\ldots,p} \, g_i(x^{(\bar{x},u,d)}(t)) > 0 \right\rbrace
    \end{align*}
    is a nonempty open subset of $[0,\infty)$.
    Therefore, we obtain that for all $u \in \mathcal{U}$, there exists a $d \in \mathcal{D}$ and a $\bar{t}(u,d) < \infty$ ($\bar{t}(u,d)$ cannot be equal to $\infty$ since otherwise $x^{(\bar{x},u,d)}$ would never leave $\cl{\RA} = \RA$, thus contradicting the assumption)  with $\max_{i = 1,2,\ldots,p} \, g_{i}(x^{(\bar{x},u,d)}(\bar{t}(u,d))) > 0$, which proves that $\bar{x} \in \RAC$, hence (ii).

    To prove (iii), by (H5), there exists $\bar{x} \in \cl{\RA}\cap \cl{\RAC} = \RA\cap \cl{\RAC} = \RAB$ such that 
    \begin{align*}
    \sup_{t \in [0,\infty)} \, \max_{i = 1,2, \ldots,p } \, g_i(x^{(\bar{x},u,d)}(t)) = 0.
    \end{align*}
    Thus, by (i) and (ii), $\bar{x}$ satisfies \eqref{eq:RA_leq_0} and \eqref{eq:RAC_greater_0} where the strict inequality is changed to $\geq$ due to the closure of $\RAC$.
   Combining these two inequalities, using (H5), we get \eqref{eq:RAB_equal_0}, which proves (iii). The converse is immediate.
\end{proof}

\begin{remark}
The  statements of Proposition \ref{Prop:Characteriztion_RAB} indeed hold true for the sets $\RAT$, $\RATC$ and $\RABT$, namely if we replace the infinite time interval $[0,\infty)$ by the finite one $[0,T]$ for all $T>0$.
\end{remark}

\begin{proposition}\label{Prop:RAB_max_equal_0}
    Assume that  assumptions (H1)-(H5) hold and that $\bar{x} \in \RAB$. Then, the supremum over $d \in \D$ in (iii) of Proposition \ref{Prop:Characteriztion_RAB} is achieved by some $\bar{d} \in \D$ and  
    \begin{align}\label{eq:RAB_max_equal_0}
    \begin{split}
		\min_{u \in \U} \, \max_{d \in \D} \, \sup_{t \in [0,\infty)} \, \max_{i = 1,2, \ldots,p } \, g_i(x^{(\bar{x},u,d)}(t)) = \sup_{t \in [0,\infty)} \, \max_{i = 1,2, \ldots,p } \, g_i(x^{(\bar{x},\bar{u},\bar{d})}(t)) = 0 .
    \end{split}
    \end{align}
\end{proposition}
\begin{proof}
In order to show that the supremum is achieved by some $\bar{d} \in \D$, a maximizing sequence $\{ d_k \}_{k\in \N}$ in $\D$ is picked such that
\begin{align*}
	\lim_{k\rightarrow \infty}  \, \sup_{t \in [0,\infty)} \, \max_{i = 1,2, \ldots,p } \, g_{i}(x^{(\bar{x},\bar{u},d_{k})}(t)) = \sup_{d \in \D} \, \sup_{t \in [0,\infty)} \, \max_{i = 1,2, \ldots,p } \, g_i(x^{(\bar{x},\bar{u},d)}(t)).
    \end{align*}
Note that $x^{(\bar{x},\bar{u},d_k)}(t)$ is entirely contained in $\RA$ for all $k \in \N$ and all $t \geq 0$, thus implying that 
$$\min_{u \in \U} \, \sup_{d \in \D} \, \sup_{t \in [0,\infty)} \, \max_{i = 1,2, \ldots,p } \, g_i(x^{(\bar{x},u,d)}(t)) \leq 0.$$
By extracting a uniformly converging subsequence, its limit, denoted by $\xi^{(\bar{u})} = \lim_{k \to \infty} x^{(\bar{x},\bar{u},d_k)}$, is absolutely continuous and admits a time derivative $\dot{\xi}^{(\bar{u})}$ in the sense of distributions. Following the same argument as in the proof of \cite[Lemma A.2]{Levine_2013},  by Mazur's Theorem (see, e.g., \cite[Chapter V, \S1, Theorem 2, p. 120]{yosida}) and the convexity assumption~(H3), there exists a convex combination of the $f(x^{(\bar{x},\bar{u},d_k)},\bar{u}, d_{k})$'s that strongly converges in $L^{2}([0,T], \R^{n})$ for all $T$ to an element of $f(\xi^{(\bar{u})}, \bar{u}, \DD)$. Therefore, as in \cite[Lemma A.2]{Levine_2013},
there exists a measurable function $\bar{d} \in \D$ such that $\xi^{(\bar{u})}(t) = x^{(\bar{x},\bar{u},\bar{d})}(t)$  and $\dot{\xi}^{(\bar{u})}(t) = f(x^{(\bar{x},\bar{u},\bar{d})}(t), \bar{u}(t),\bar{d}(t))$  for almost all $t \geq 0$.
Moreover we have:
\begin{align*}
\lim_{k \to \infty} \sup_{t \in [0,\infty)} \max_{i = 1,2, \ldots,p } g_i\big(x^{(\bar{x},\bar{u},d_k)}(t)\big) 
= \sup_{t \in [0,\infty)} \max_{i = 1,2, \ldots,p } g_i\big(x^{(\bar{x},\bar{u},\bar{d})}(t)\big) = \sup_{d \in \D} \, \sup_{t \in [0,\infty)} \, \max_{i = 1,2, \ldots,p } \, g_i(x^{(\bar{x},\bar{u},d)}(t)) = 0.
\end{align*}
which proves \eqref{eq:RAB_max_equal_0} and that the supremum over $\D$ is achieved by $\bar{d}$, hence the Proposition.
\end{proof}

As in \cite{Levine_2013}, we  introduce the following splitting of the boundary $\RAB$ into two disjoint subsets, the  \emph{usable part of the boundary}, or more simply, the \emph{usable part}, and the  \emph{barrier}, a terminology borrowed from \cite{Isaacs_1965}.

\begin{definition}
    The set $\RAO \triangleq \RAB \cap \G_0$ is called the \emph{usable part} of $\RAB$ and the set $\RAM \triangleq \RAB \cap \G_-$ is called the \emph{barrier} of $\RA$.
\end{definition}

Since $\RA \subset \G = \G_0 \cup \G_-$, we immediately obtain $\RAB = \RAO  \cup \RAM$.

\subsection{Tangent Property of the Usable Part}

To further characterize $\RAO$, we introduce the Lie-derivative of a $\mathcal{C}^1$ function $h: \R^n \rightarrow \R$ along the vector field $f(\cdot,u,d)$ at a point $x \in \R^n$, denoted by $L_f h(x,u,d) \triangleq D h(x) f(x,u,d)$, where $D h(x)$ is the gradient of $h$ at the point $x$.

We next introduce the set of active indices at a point $x \in \R^n$, \ie the set of indices  $i\in \{1, \ldots, p\}$ for which the state constraint is active, namely $g_i(x)=0$:
\begin{definition}
    At a given point $x \in \R^n$, the set of active indices, denoted by  $\II(x)$, of the state constraints $g$ is given by
    \begin{align*}
        \II(x) \triangleq \left\lbrace i \in \{1,2,\ldots,p \} \; \vert \; g_i(x) = 0 \right\rbrace . 
    \end{align*}
\end{definition}
Indeed, if $x\not\in \G_0$, $ \II(x) = \emptyset$.

The usable part $\RAO$ enjoys the following tangent property:
\begin{proposition}\label{Prop:Characterizing_RAO}
    $\RAO$ is contained in the set of points $z \in \G_0$ such that
    \begin{align*}
        \min_{u \in \U} \, \max_{d \in \D}  \, L_f g_i(z,u,d) \leq 0, \quad  \forall i \in \II(z).
    \end{align*}
\end{proposition}
 \begin{proof}
Consider a point $z \in \RAO$.
    Thus, there exists a control $u \in \mathcal{U}$ such that for all disturbances $d \in \D$ the trajectory $x^{(x_0,u,d)}$ remains in $\RA$ for all time.
    In particular, at time $t = 0$, according to Proposition~\ref{Prop:RAB_max_equal_0}, the right limit of the vector $f(z,u(t),d(t))$, \ie $\lim_{t \downarrow 0} \, f(z,u(t),d(t))$, cannot point outwards of $\mathcal{G}$ for all $d \in \mathcal{D}$.
    Otherwise, since $g_{i}$ is $C^2$ by (H4), if there would exist an $i \in \mathbb{I}(z)$ and a $d\in \D$ such that for all $u\in \U$,  $\lim_{t \downarrow 0} \frac{d}{dt}g_{i}(x^{(z,u,d)}(t)) = \lim_{t \downarrow 0} \, L_f g_{i}(z,u(t),d(t)) > 0$, the resulting trajectory would leave the constraint set $\mathcal{G}$ at least in a small interval of time which contradicts the fact that $z \in \RAO$.
    Therefore, for all $i \in \mathbb{I}(z)$ and all $d \in \mathcal{D}$ we must have $\lim_{t \downarrow 0} \, L_f g_i(z,u(t),d(t)) \leq 0$.
    This implies $\lim_{t \downarrow 0} \, \sup_{d \in \mathcal{D}} \, L_f g(z,u(t),d(t)) \leq 0$ for all $i \in \mathbb{I}(z)$, which proves the statement.
\end{proof}

\subsection{Semi-Permeability of the Barrier}
Let us now focus on the barrier $\RAM$. We adapt to our context the  dynamic programming argument of the proof of Proposition 5.3 in \cite{Levine_2013}, in order to prove that the barrier is comprised of segments of system's integral curves until their intersection with $\G_0$.
\begin{proposition}\label{Prop:Characterizing_Trajectories_in_RAM}
    Assume that assumptions (H1)-(H5) hold.
    The set $\RAM$ is made of points $\bar{x} \in \G_-$ for which there exists a control input $\bar{u} \in \U$ and a disturbance $\bar{d} \in \D$, such that the resulting integral curve $x^{(\bar{x},\bar{u}, \bar{d})}$ is entirely contained in $\RAM$ until it intersects with $\G_0$ at $x^{(\bar{x},\bar{u}, \bar{d})}(\bar{t})$ for some $\bar{t} \in [0, \infty)$.
\end{proposition}
\begin{proof}
Assume that  $\bar{x} \in \RAM$.
    Thus, $\bar{x} \in \RAB$ and, by Proposition \ref{Prop:RAB_max_equal_0}, there exists a $\bar{u} \in \mathcal{U}$ and a $\bar{d} \in \mathcal{D}$ such that $\sup_{t \in [0,\infty)} \, \max_{i = 1,2, \ldots,p } \, g_i(x^{(\bar{x},\bar{u},\bar{d})}(t)) = 0$.
    Moreover, there exists a $\bar{t} \in [0, \infty)$ such that the resulting integral curve $x^{(\bar{x},\bar{u},\bar{d})}$ is in $\mathcal{G}_-$ for all $t \in [0,\bar{t})$ and reaches $\mathcal{G}_0$ at $\bar{t}$ for the first time, \ie, $\max_{i = 1,2, \ldots,p } \, g_i(x^{(\bar{x},\bar{u},\bar{d})}(t)) < 0$ for all $t \in [0,\bar{t})$ and $\max_{i = 1,2, \ldots,p } \, g_i(x^{(\bar{x},\bar{u},\bar{d})}(\bar{t})) = 0$.
    Now consider any $t_0 \in [0,\bar{t})$ and let $\xi = x^{(\bar{x},\bar{u},\bar{d})}(t_0)$.
    Note that $\xi \in \mathcal{G}_-$.
    With a standard dynamic programming argument and using $t_0 < \bar{t}$, we obtain
$$
\min_{u \in \U} \, \max_{d \in \D} \, \sup_{t \in [0,\infty)} \, \max_{i = 1,2, \ldots,p } \, g_i(x^{(\xi,u,d)}(t)) =  
\min_{u \in \U} \, \max_{d \in \D} \, \sup_{t \in [0,\infty)} \, \max_{i = 1,2, \ldots,p } \, g_i(x^{(\bar{x},u,d)}(t+t_0)) = 0
$$
    which implies that $\xi \in \RAM$ by Proposition \ref{Prop:RAB_max_equal_0}.
    Since $t_0 \in [0,\bar{t})$ was arbitrary, the whole arc of the integral curve starting from $\bar{x} \in \RAM$ is contained in $\RAM$ until it intersects with $\mathcal{G}_0$ at the time $\bar{t}$, which proves the proposition.
 \end{proof}

For the next results, we need to introduce the concatenation of two controls.
For $u_1, u_2 \in \U$ and $\tau \geq t_0$, the concatenated control input at $\tau$, denoted by $v \triangleq  u_1 \Join_\tau u_2$, is given by
\begin{align*}
    v(t) \triangleq u_1 \Join_\tau u_2 =  \begin{cases}
    u_1(t) & \mathrm{if~} t \in [t_0, \tau) \\
    u_2(t) & \mathrm{if~} t \geq \tau .
    \end{cases}
\end{align*}
Note that $v \in \U$.
Moreover, it is readily seen that every integral curve $x^{(x_0,t_0,v,d)}$ of system \eqref{eq:control_system} with $d \in \D$ and $\xi = x^{(x_0,t_0,u_1,d)}(\tau)$ satisfies
\begin{align*}
    x^{(x_0,t_0,v,d)} = x^{(x_0,t_0,u_1,d)} \Join_\tau x^{(\xi,\tau,u_2,d)},
\end{align*}
where, with a slight abuse of notation, $\Join_\tau$ also refers to the concatenation of two integral curves at time $\tau$.

The same concatenation definition also holds for disturbances $d$ and pairs $(u,d)$ as well as their respective integral curves.

\begin{corollary}\label{Cor:Trajectories_remain_in_RAM}
    From any point of the barrier $\bar{x} \in \RAM$, there exists $\bar{d}\in \D$ such that no $u \in \U$ can force the integral curve $x^{(\bar{x},u,\bar{d})}$ to penetrate $\textnormal{int}(\RA)$ before leaving $\G_-$.
\end{corollary}
\begin{proof}
 The proof is by contradiction.
    Consider a point $\bar{x} \in \RAM$.
    From Proposition \ref{Prop:Characterizing_Trajectories_in_RAM}, there exists a control input $\bar{u} \in \mathcal{U}$ together with a maximizing disturbance $\bar{d} \in \mathcal{D}$ such that $x^{(\bar{x},\bar{u},\bar{d})}(t) \in \RAM$ for all $t \in [0,\bar{t})$, where $\bar{t} \in [0,\infty)$ is the first time such that $x^{(\bar{x},\bar{u},\bar{d})}(\bar{t}) \in \mathcal{G}_0$.
    Let $\tau \in [0,\bar{t})$ and $\xi = x^{(\bar{x},\bar{u},\bar{d})}(\tau)$.
    Moreover, let $\varepsilon >0$ small enough such that $x^{(\bar{x},\bar{u},\bar{d})}(\tau + t) \in \mathcal{G}_-$ for all $t \in [0,\varepsilon]$.
    Assume that there exists a control $\tilde{u} \in \mathcal{U}$ and a time $\sigma \in [0,\varepsilon)$ with $\zeta = x^{(\xi,\tilde{u},d)}(\tau + \sigma) \in \textnormal{int}(\RA)$ for all $d \in \mathcal{D}$.
    Thus, by Proposition \ref{Prop:Characteriztion_RAB} (i) and (iii) as well as Proposition \ref{Prop:RAB_max_equal_0}, we obtain that there exists a control $\tilde{v} \in \mathcal{U}$ such that 
    \begin{align*}
        \sup_{d \in \mathcal{D}} \, \sup_{t \in [\tau + \sigma, \infty)} \, \max_{i = 1,2,\ldots,p} \, g(x^{(\zeta,\tilde{v},d)}(t)) < 0.
    \end{align*}
    By concatenating $\bar{u}$, $\tilde{u}$ and $\tilde{v}$ as $v:= \bar{u} \Join_\tau \tilde{u} \Join_{\tau + \sigma} \tilde{v}$, we immediately obtain that
$$
\sup_{d \in \D} \, \sup_{t \in [0, \infty)} \, \max_{i = 1,2,\ldots,p} \, g_{i}(x^{(\bar{x},v,d)}(t)) = 
\sup_{d \in \D} \, \sup_{t \in [\tau + \sigma, \infty)} \, \max_{i = 1,2,\ldots,p} \, g_{i}(x^{(\zeta,v,d)}(t)) < 0,
$$
    but since $\bar{x} \in \RAM \subset \RAB$, the latter inequality contradicts Proposition \ref{Prop:Characteriztion_RAB} (iii).
\end{proof}

\begin{remark}
    Note that Proposition \ref{Prop:Characterizing_Trajectories_in_RAM} does not say more on the behavior of the integral curve $x^{(\bar{x},\bar{u},\bar{d})}$, after its intersection with $\G_0$ at the time $\bar{t}$, than the fact that it remains in $\G$.
    It could possibly remain in $\G_0$ for some time until it returns to $\G_-$, or could also remain in $\G_0$ forever and  get stuck.
\end{remark}

\begin{remark}
    Corollary \ref{Cor:Trajectories_remain_in_RAM} may be interpreted as a \emph{semi-permeability property} according to Isaacs's terminology \cite{Isaacs_1965}: starting in  $\RAM$, one can leave $\RA$ and $\G_-$ after some time, but the worst-case disturbances prevent the state from going back into $\RA$. Hence, for such disturbances, $\RAM$ can be crossed only in one direction, namely towards the complement of $\G_-$.
\end{remark}

\section{The Barrier's Ultimate Locally Separating Hyperplane}\label{Sec:Ultimate_Sep}

We now study the local tangential structure of the barrier $\RAM$ at its intersection with the active state constraints $\G_0$. We show the necessity of the existence of a separating hyperplane at this intersection that will serve later on as a boundary condition for a min-max version of the Pontryagin maximum  principle.

We assume from now on that the set $\LL_0$ defined by
\begin{equation}\label{LL0set:def}
\LL_0 \triangleq \{ (z,u,d)\in \G_0 \times \UU\times \DD \mid  \, L_f g_{i}(z,u,d) = 0, \; i \in \II(z)\}
\end{equation}
is nonempty. 

$\LL_0$ may be interpreted as the subset of $\G_0 \times \UU\times \DD$ where 
the vector  $f(z,u,d)$ is tangent to $\G_0$.

More precisely, if $(z,\bar{u},\bar{d})\in \LL_0$, then $z$ and $i$ satisfy: 
\begin{equation}\label{LL0set:eq}\max_{i=1,\ldots, p} g_{i}(z)=0 , \quad   L_f g_{i}(z,\bar{u},\bar{d}) = 0, \quad i \in \II(z),
\end{equation}
where none of $i$, $z$, $\bar{u}$, $\bar{d}$ is necessarily unique.

Note that the case $n=1$ is trivial since it reads $\dot{x}=v$, subject to $g_{i}(x)\leq 0$, $i=1,\ldots, p$, after re-mapping $u$ and $d$ to $v=f(x,u,d) \in \R$ and $L_f g_{i}(z, u, d) = \frac{dg_{i}}{dx}(z) v =0$ for  $i\in \II(z)$ and $v\in f(z, \UU,\DD)$, \ie $\frac{dg_{i}}{dx}(z)=0$ or $v=0$ if $0\in f(z, \UU,\DD)$.

Therefore, in the sequel, we assume that $n\geq 2$ and that the projection of $\LL_0$:
$$\pi \LL_0 \triangleq \{ z\in \G_0 \mid \exists (\bar{u},\bar{d})\in \UU\times \DD : (z, \bar{u},\bar{d})\in \LL_0\},$$ 
\ie satisfying \eqref{LL0set:eq}. We further assume the following
\begin{itemize}
\item[(H6)] $\pi \LL_0$ is locally $C^1$-diffeomorphic to an $(n-2)$-dimensional polytope.
\end{itemize}

\begin{remark}\label{rem:H6}
Assumption (H6) guarantees that the equation $L_f g_{i}(z,\bar{u},\bar{d}) = 0$, with $i \in \II(z)$, is well-defined, as well as the sets $\LL_0$ and $\pi \LL_0$:
indeed, if $i,j\in \{1, \ldots, p\}$, $i\neq j$ and $z\in \G_0$, are such that $g_{i}(z)= g_{j}(z) = \max_{l= 1, \ldots, p} g_{l}(z) = 0$, then $g_{i}$ and $g_{j}$ necessarily have colinear tangent vectors at this point and therefore $L_{f}g_{i}(z, u, d)= L_{f}g_{j}(z, u, d)$. Therefore the definition~\eqref{LL0set:eq} makes sense even if the number of elements of $\II(z)$ is larger than 1.
\end{remark}

\subsection{The Ultimate Locally Separating Hyperplane Condition}

At the intersection of a trajectory in $\RAM$ with $\G_0$, we now prove that the control and disturbance must satisfy a local min–max (saddle-point) condition defining the ultimate locally separating hyperplane of the barrier, separating the robust admissible set $\RA$ and its complement $\RAC$ at this intersection point.
\begin{proposition}\label{Prop:Ultimate_Separation}
We assume that (H1)--(H6) hold true and we consider $\bar{x} \in \RAM$ as well as a minimizing control law $\bar{u} \in \U$ and a maximizing disturbance $\bar{d} \in \D$ satisfying Proposition~\ref{Prop:Characterizing_Trajectories_in_RAM}, \ie  the integral curve $x^{(\bar{x},\bar{u}, \bar{d})}$ runs along $\RAM$ and intersect with $\G_0$ at some time $\bar{t} \in [0,\infty)$.
    Then the intersection point $z = x^{(\bar{x},\bar{u}, \bar{d})}(\bar{t})$ is such that, up to a suitable change of $\bar{u}$ and $\bar{d}$ on a 0-measure set,
        \begin{equation}\label{LL0z:eq}
        (z,\bar{u}(\bar{t}),\bar{d}(\bar{t}))\in \LL_0
        \end{equation} 
        and $(\bar{u}(\bar{t}),\bar{d}(\bar{t}))$ is a saddle-point of the mapping 
        $$(u,d)\in \UU \times \DD \mapsto  \, L_f g_{i}(z,u,d), $$ 
        i.e.
        \begin{align}\label{eq:saddle}
        \begin{split}
            L_f g_{i}(z,\bar{u}(\bar{t}),d) &\leq \, L_f g_{i}(z,\bar{u}(\bar{t}),\bar{d}(\bar{t})) =0 \leq \, L_f g_{i}(z,u,\bar{d}(\bar{t})), \quad \forall i \in \II(z), \quad \forall (u,d) \in \UU \times \DD.
        \end{split}
        \end{align}
\end{proposition}

\begin{proof}
Consider, as in Proposition \ref{Prop:Characterizing_Trajectories_in_RAM},, an initial point $\bar{x} \in \RAM$, a minimizing control law $\bar{u} \in \U$ and a maximizing disturbance $\bar{d} \in \D$ such that the integral curve runs along $\RAM$ and intersects with $\G_0$ at some time $\bar{t} \in [0,\infty)$. Two, and only two, cases are indeed possible: $\bar{t} < \infty$ or $\bar{t} = \infty$.
    
\subsubsection{The case $\bar{t}<\infty$} Let us first assume that $\bar{t}$ is finite and set $z = x^{(\bar{x},\bar{u}, \bar{d})}(\bar{t})$.     
    
According to \eqref{eq:RAB_max_equal_0}, we have $\max_{i \in \{1, \ldots,p\}} \, g_{i}(x^{(\bar{x},\bar{u}, \bar{d})}(\bar{t})) = 0= \max_{i\in \II(z)}g_{i}(z)$ and, according to \eqref{eq:RA_leq_0}, we have
\begin{align*}
    \max_{i \in \{1, \ldots,p\}} \, g_{i}(x^{(\bar{x},\bar{u}, d)}(\bar{t})) \leq \sup_{t \in [0, \infty)} \, \max_{i \in \{1, \ldots,p\}} \, g_{i}(x^{(\bar{x},\bar{u}, d)}(t)) &\leq \max_{d\in \D}\sup_{t \in [0, \infty)} \, \max_{i \in \{1, \ldots,p\}} \, g_{i}(x^{(\bar{x},\bar{u}, d)}(t)) \\
    &\leq \max_{i \in \{1, \ldots,p\}} \, g_{i}(x^{(\bar{x},\bar{u}, \bar{d})}(\bar{t})) = 0
\end{align*}
for all $d \in \D$.
Now, letting $\tilde{d}$ be an arbitrary needle perturbation of $\bar{d}$ at the Lebesgue perturbation time $\bar{t}-l\varepsilon$ (see Appendix \ref{Appendix:Needle_Perturbations}): 
$$\tilde{d} = \bar{d} \Join_{\bar{t}-l\varepsilon} d \Join_{\bar{t}}  \bar{d}.$$ 
According to \eqref{eq:Lemma_Needle_Perturb_Approx}, we get, after division by $\varepsilon$,
\begin{align*}
    \max_{d \in \DD}  \, D g_{i} (z) f(z, \bar{u}(\bar{t}), d) + 0(\varepsilon) = \max_{d \in \DD} \, \frac{1}{\varepsilon}\left( g_{i}(x^{(\bar{x},\bar{u}, \tilde{d})}(\bar{t})) - g_{i}(z) \right) \leq 0, \quad i \in \II(z)
\end{align*}
and, taking the limit with respect to $\varepsilon$, we conclude that, for $i \in \II(z)$,
\begin{align}\label{eq:max_d_Lfg_leq_0} 
    \max_{d \in \DD} \, L_f g_{i}(z,\bar{u},d) \leq  \, L_f g_{i}(z,\bar{u},\bar{d}) \leq 0.
\end{align}

Symmetrically, interchanging the roles of $d$ and $u$, and taking $\tilde{u}$ as an arbitrary needle perturbation of $\bar{u}$ at the Lebesgue perturbation time $\bar{t}-l\varepsilon$, \ie $\tilde{u} = \bar{u} \Join_{\bar{t}-l\varepsilon} u \Join_{\bar{t}}  \bar{u}$, according to \eqref{eq:Lemma_Needle_Perturb_Approx}, we get
\begin{align*}
    0 = \max_{i \in \{1, \ldots,p\}} \, g_{i}(x^{(\bar{x},\bar{u}, \bar{d})}(\bar{t})) = \min_{u \in \U} \sup_{t\in[0, \infty)} \max_{i \in \{1, \ldots,p\}} \, g_{i}(x^{(\bar{x},u, \bar{d})}(t)) \leq \sup_{t\in[0, \infty)} \max_{i \in \{1, \ldots,p\}} \, g_{i}(x^{(\bar{x},u, \bar{d})}(t))
\end{align*}
or
\begin{align*}
    0 \leq \min_{u \in \UU} \, \max_{i \in \II(z)} \, \frac{1}{\varepsilon}\left( g_{i}(x^{(\bar{x},\tilde{u}, \bar{d})}(\bar{t})) - g_{i}(z) \right) = \min_{u \in \UU} \, \max_{i \in \II(z)} \, Dg_{i} (z) f((z, u, \bar{d}(\bar{t})) + 0(\varepsilon) = \min_{u \in \UU} \, \max_{i \in \II(z)} \, L_{f} g_{i} (z, u, \bar{d}(\bar{t})) + 0(\varepsilon)
\end{align*}
for all $\ee$ sufficiently small and for all $u\in \UU$ (recall that $\tilde{u} = \bar{u} \Join_{\bar{t}-l\varepsilon} u \Join_{\bar{t}}  \bar{u}$). Thus, for $i \in \II(z)$,
\begin{align}\label{eq:0_leq_min_u_Lfg} 
    0 \leq \min_{u \in \UU} \,  L_{f} g_{i} (z, u, \bar{d}(\bar{t})).
\end{align}
Hence, grouping \eqref{eq:max_d_Lfg_leq_0} and \eqref{eq:0_leq_min_u_Lfg}, we have proven the saddle-point property \eqref{eq:saddle} with \eqref{LL0z:eq}.

\bigskip
\subsubsection{The case $\bar{t}=\infty$}
Now, assume that $\bar{t}$ is infinite and set $z = \lim_{t \rightarrow \infty} \, x^{(\bar{x},\bar{u}, \bar{d})}(t)$ with $x^{(\bar{x},\bar{u}, \bar{d})}(t)\in \RAM$ for all $t\geq 0$.
If $\{t_k\}_{k\in \N}$ is an arbitrary increasing sequence of time instants made of Lebesgue points of $\bar{u}$ and $\bar{d}$, since $\bar{u}(t_{k})$ and $\bar{d}(t_{k})$ are bounded and $\UU\times\DD$ is compact, there exists a  subsequence, still noted $(\bar{u}(t_{k}),\bar{d}(t_{k}))$, converging to some $(\hat{u}, \hat{d})\in \UU\times \DD$. Thus, by the continuity of $x\mapsto \max_{i \in \{1, \ldots,p\}} g_{i}(x)$, we get 
$$\max_{i \in \{1, \ldots,p\}} g_{i}(z) = \lim_{t_{k}\rightarrow \infty} \max_{i \in \{1, \ldots,p\}} g_{i}(x^{(\bar{x},\bar{u}, \bar{d})}(t_{k})) =\sup_{t \in [0, \infty)} \, \max_{i \in \{1, \ldots,p\}} \, g_{i}(x^{(\bar{x},\bar{u}, \bar{d})}(t)) =0$$
and for all $\ee$ small enough and $k$ large enough, we have
\begin{align*}
\begin{split}
\Big\vert \max_{i \in \{1, \ldots,p\}} g_{i}(z) - \max_{i \in \{1, \ldots,p\}} g_{i}(x^{(\bar{x},\bar{u}, \bar{d})}(t_{k})) \Big\vert \leq \Big\vert \max_{i \in \{1, \ldots,p\}} Dg_{i}(z)(z-x^{(\bar{x},\bar{u}, \bar{d})}(t_k)) \Big\vert \leq \ee.
\end{split}
\end{align*}

 Hence, for $i\in \II(z)$,
 \begin{equation}\label{lim-L0:eq}
 \lim_{k\rightarrow \infty} Dg_{i}(z)(z-x^{(\bar{x},\bar{u}, \bar{d})}(t_k)) = Dg_{i}(z)f(z,\hat{u}, \hat{d})=0
 \end{equation}
which proves that  
$(z,\hat{u},\hat{d})\in \LL_0$.

We now prove the saddle point property \eqref{eq:saddle}. We have
\begin{align*}
\min_{u\in \U} \max_{d\in \D} \sup_{t \in [0, \infty)} \, \max_{i \in \{1, \ldots,p\}} \, g_{i}(x^{(\bar{x},u,d)}(t)) = \sup_{t \in [0, \infty)} \, \max_{i \in \{1, \ldots,p\}} \, g_{i}(x^{(\bar{x},\bar{u}, \bar{d})}(t)) = 0.
\end{align*}

We again let $\tilde{d}$ be an arbitrary needle perturbation of $\bar{d}$ at the perturbation time $t_{k} - l \ee$, \ie $\tilde{d} = \bar{d} \Join_{t_{k}-l\ee} d \Join_{t_{k}}  \bar{d}$ for $t_{k}$ large enough.
Thanks to Lemma \eqref{eq:Lemma_Needle_Perturb_Approx}, for all $t\in ]t_{k}-l\ee ,t_{k}[$ and all $d\in \DD$, we have
\begin{align*}
  x^{(\bar{x},\bar{u}, \tilde{d})}(t) - x^{(\bar{x},\bar{u}, \bar{d})}(t) = \ee \left( f(x^{(\bar{x},\bar{u}, \tilde{d})}(t_k) - f(x^{(\bar{x},\bar{u}, \bar{d})}(t_k) \right) +0(\ee^2)
\end{align*}
Therefore we have, again for all $t\in ]t_{k}-l\ee ,t_{k}[$ with $k$ large enough, all $d\in \DD$ and $i\in \II(x^{(\bar{x},\bar{u}, \bar{d})}(t))$,
\begin{align*}
g_{i}(x^{(\bar{x},\bar{u}, \tilde{d})}(t)) - g_{i}(x^{(\bar{x},\bar{u}, \bar{d})}(t)) \ &=  Dg_{i}(x^{(\bar{x},\bar{u}, \bar{d})}(t))\left (x^{(\bar{x},\bar{u}, \tilde{d})}(t) - x^{(\bar{x},\bar{u}, \bar{d})}(t)\right) +0(\ee)\\
  &=\ee Dg_{i}(x^{(\bar{x},\bar{u}, \bar{d})}(t))\left ( f(x^{(\bar{x},\bar{u}, \tilde{d})}(t_k),\bar{u}(t_k), \tilde{d}(t_k)) - f(x^{(\bar{x},\bar{u}, \bar{d})}(t_k) , \bar{u}(t_{k}), \bar{d}(t_{k}) )\right) + 0(\ee^2)\\
  &=\ee Dg_{i}(x^{(\bar{x},\bar{u}, \bar{d})}(t)) f(x^{(\bar{x},\bar{u}, \tilde{d})}(t_k),\bar{u}(t_k), \tilde{d}(t_k)) + 0(\ee^2)
    \end{align*} 
    since $Dg_{i}(x^{(\bar{x},\bar{u}, \bar{d})}(t))f(x^{(\bar{x},\bar{u}, \bar{d})}(t_k) , \bar{u}(t_{k}), \bar{d}(t_{k}) ) = 0(\ee^2)$ by \eqref{lim-L0:eq} and thus 
$$
  g_{i}(x^{(\bar{x},\bar{u}, \tilde{d})}(t)) - g_{i}(x^{(\bar{x},\bar{u}, \bar{d})}(t)) = \ee L_f g_{i}(x^{(\bar{x},\bar{u}, \bar{d})}(t),\bar{u},d) + 0(\ee^2).
$$

Thus, since
\begin{align*}
  \max_{d\in \DD} \max_{i \in \{1, \ldots,p\}} \, \frac{1}{\ee}\left( g_{i}(x^{(\bar{x},\bar{u}, \tilde{d})}(t)) - g_{i}(x^{(\bar{x},\bar{u}, \bar{d})}(t)) \right) \leq 0,
\end{align*}

taking the limit as $t\to \infty$ to infinity (\ie as $k\to \infty$), and as $\varepsilon \to 0$, one obtains, for $i \in \II(z)$,
\begin{align}
\begin{split}\label{eq:inf_t_max_d_Lfg_leq_0}
    \max_{d \in \DD}  L_f g_{i}(z,\hat{u},d) = \lim_{\varepsilon \rightarrow 0} \lim_{t \rightarrow \infty} \max_{d \in \DD} D g_{i} (x^{(\bar{x},\bar{u}, \tilde{d})}(t)) f(x^{(\bar{x},\bar{u}, \tilde{d})}(t), \bar{u}(t), d) + \varepsilon \leq L_f g_i(z,\hat{u},\hat{d}) = 0.
\end{split}
\end{align}
Symmetrically, interchanging the roles of $d$ and $u$, and taking $\tilde{u}$ as an arbitrary needle perturbation of $\bar{u}$ at the perturbation time $\tilde{t} - \varepsilon l$, \ie $\tilde{u} = \bar{u} \Join_{\tilde{t} - l \varepsilon} u \Join_{\tilde{t}}  \bar{u}$, following the same lines as above, we obtain, for $i \in \II(z)$,
\begin{align}
\begin{split}\label{eq:inf_t_0_leq_min_u_Lfg}
    \min_{u \in \UU} \, L_{f} g_{i} (z), u, \hat{d}) = \lim_{\varepsilon \rightarrow 0} \, \lim_{t \rightarrow \infty} \, \min_{u \in \UU} \, L_{f} g_{i} (x^{(\bar{x},\bar{u}, \bar{d})}(t), u, \bar{d}(t)) + \varepsilon \geq \, L_{f} g_{i} (z, \hat{u}, \hat{d}) =0.
\end{split}
\end{align}
By combining \eqref{eq:inf_t_max_d_Lfg_leq_0} and \eqref{eq:inf_t_0_leq_min_u_Lfg}, we have shown the saddle-point property \eqref{eq:saddle} when $\bar{t}= +\infty$.
\end{proof}
\begin{remark}
 The limits of $\bar{u}(t)$ and $ \bar{d}(t)$  do not necessarily exist but $\hat{u}$ and $\hat{d}$ are obtained as in the proof of Proposition \ref{Prop:RAB_max_equal_0}, as accumulation points of convex combinations of the elements of subsequences of $f(z,\bar{u}(t_k),\bar{d}(t_k))\in f(z, \UU,\DD)$ which is convex by (H3).
\end{remark}
\begin{remark}
    The saddle point property \eqref{eq:saddle} at the intersection point $z$ of an integral curve on the boundary of the admissible set with the constraints has the following geometric interpretation: the linear form $f\mapsto L_f g_{i}(z,\bar{u},\bar{d}) =\langle Dg_{i}(z),f(z,\bar{u},\bar{d})\rangle$, with $i\in \II(z)$, is the separating hyperplane between the admissible set $(\max_{u\in \UU}\langle Dg_{i}(z),f(z,\bar{u},\bar{d})\rangle\leq 0)$  and its complement $(\min_{d\in \DD}\langle Dg_{i}(z),f(z,\bar{u},\bar{d})\rangle\geq 0)$ and is tangent to the barrier at the point $z$ ($\langle Dg_{i}(z),f(z,\bar{u},\bar{d})\rangle =0$).
    \end{remark}

\begin{remark}
Note that the condition \eqref{eq:saddle} slightly differs from standard existence theorems of saddle points as can be seen, \eg, in \cite{Von_Neumann_1944,Isaacs_1965,Friedman_1970,Berkovitz_1985,Krassovski_1977}, where it is assumed that the functional $(u,d) \mapsto \langle Dg_{i}(z),f(z,u,d)\rangle$  satisfies
$u \mapsto \langle Dg_{i}(z),f(z,u,d)\rangle$: convex and lower semi-continuous on the compact $\UU$ for all $d$ and
$d \mapsto \langle Dg_{i}(z),f(z,u,d)\rangle$: concave and upper semi-continuous on the compact $\DD$ for all $u$. Here, this functional is assumed to be continuous in both variables but the convexity comes from the assumption (H3) on the images of $\UU$ and $\DD$ by $f$.
\end{remark}

\subsection{The barrier equation}
The following theorem is a generalization in the robust admissible set's context of the main result of \cite{Levine_2013}[Theorem 7.1]: the barrier must be made of integral curves satisfying a saddle point version of Pontryagin's  maximum principle, where the final condition of the adjoint is given by the ultimate separating hyperplane condition.

\begin{theorem}\label{Th:Barrier_Theorem}
    Under the assumptions (H1)-(H5), every triple $(\bar{u}, \bar{d}, x^{(\bar{u},\bar{d})}) \in \U \times \D \times \left(\RAM \cap \cl{\Int{\RA}}\right)$ satisfy the following necessary conditions:

    There exists a nonzero absolutely continuous solution $\lambda^{(\bar{u},\bar{d})}$ to the adjoint equation
    \begin{align}\label{eq:Th_adjoint_system}
    \begin{split}
        \dot{\lambda}^{(\bar{u},\bar{d})}(t) &= - \left( \frac{\partial f}{\partial x} (x^{(\bar{u},\bar{d})}(t), \bar{u}(t), \bar{d}(t) ) \right)^\top \lambda^{(\bar{u},\bar{d})}(t), \\
        \lambda^{(\bar{u},\bar{d})}(\bar{t}) &= \left(D g_{i}(z) \right)^\top, \quad i \in \mathbb{I}(z),
    \end{split} 
    \end{align}
such that
    \begin{align}
    \begin{split}\label{eq:Th_Hamiltonian}
        \min_{u \in \UU} \, \max_{d \in \DD} \, \left\lbrace (  \lambda^{(\bar{u},\bar{d})}(t))^\top f(x^{(\bar{u},\bar{d})}(t), u,d) \right\rbrace &= \max_{d \in \DD} \, \min_{u \in \UU} \, \left\lbrace (\lambda^{(\bar{u},\bar{d})}(t))^\top f(x^{(\bar{u},\bar{d})}(t), u,d) \right\rbrace \\
        &= (\lambda^{(\bar{u},\bar{d})}(t))^\top f(x^{(\bar{u},\bar{d})}(t), \bar{u}(t), \bar{d}(t)) = 0
    \end{split}
    \end{align}
    at every Lebesgue point $t$ of $(\bar{u},\bar{d})$, \ie, for almost all $t \leq \bar{t}$,
    with $z$ satisfying $z = x^{(\bar{u},\bar{d})}(\bar{t}) \in \G_0$ given by
    \begin{align}\label{eq:Th_UT-Condition}
    g_i(z) = 0, \quad
    \min_{u \in \UU} \max_{d \in \DD}  L_f g_i(z,u,d) &= \max_{d \in \DD} \min_{u \in \UU}  L_f g_i(z,u,d) = L_f g_{i}(z,\bar{u},\bar{d}) = 0, \quad \forall i \in \II(z).
    \end{align}
    Moreover, $\lambda^{(\bar{u},\bar{d})}(t)$ is normal to $\RAM \cap \textnormal{cl}(\textnormal{int}(\RA))$ at $x^{(\bar{u},\bar{d})}(t)$ for almost every $t \leq \bar{t}$, $\bar{t}$ being arbitrary.
\end{theorem}

\begin{proof}
Let us first fix $\bar{d}$ and consider the reachable set $X_{t}(\bar{x},\bar{d}) \triangleq \{ x^{(\bar{x},u,\bar{d})}(t) \, \vert \, u \in \U \}$.
We know from \cite{Levine_2013}[Proposition 7.1] that points on $\RAM \cap \textnormal{cl}(\textnormal{int}(\RA))$ are also on $\partial X_{t}(\bar{x},\bar{d})$ for some $t$ and therefore that there exists an adjoint $\eta^{\bar{d}}$ satisfying (see \cite{Levine_2013}[Theorem 7.1])
\begin{align}\label{eq:Th_adjoint-d-_system}
    \begin{split}
\dot{\eta}^{\bar{d}}(t) &= - \left( \frac{\partial f}{\partial x} (x^{(\bar{u},\bar{d})}(t), \bar{u}(t), \bar{d}(t) ) \right)^\top \eta^{\bar{d}}(t), \\  \eta^{\bar{d}}(\bar{t}) &= \left(D g_{i^*}(z) \right)^\top
    \end{split}
    \end{align}
    and
    \begin{align}\label{eq:Th_Hamilton-d-_system}
    \begin{split}
\min_{u\in \UU}  \eta^{\bar{d}}(t)^\top f(x^{(\bar{u},\bar{d})}(t), u, \bar{d}(t) ) = \eta^{\bar{d}}(t)^\top f(x^{(\bar{u},\bar{d})}(t), \bar{u}(t), \bar{d}(t) ) =0
 \end{split}
    \end{align}
for almost all $t\leq \bar{t}$, 
where $\bar{t}$ denotes the time at which $z = x^{(\bar{u},\bar{d})}(\bar{t})$ is reached, with $z \in \G_0$ satisfying
    \begin{align}\label{eq:Th_UT-d-Condition}
    \begin{split}
    g_i(z) &= 0, \; \; i \in \II(z), \\
    \min_{u \in \UU}  \max_{i \in \II(z)} L_f g_i(z,u,\bar{d}) 
	&= L_f g_{i^*}(z,\bar{u},\bar{d}) = 0.
     \end{split}
    \end{align}

Now, fixing $\bar{u}$ and considering the reachable set $X_{t}(\bar{x},\bar{u}) \triangleq \{ x^{(\bar{x},\bar{u},d)}(t) \, \vert \, d \in \D \}$ the points on $\RAM \cap \textnormal{cl}(\textnormal{int}(\RA))$ are also on $\partial X_{t}(\bar{x},\bar{u})$ for some $t$ and therefore there exists an adjoint $\eta^{\bar{u}}$ satisfying (see again \cite{Levine_2013}[Theorem 7.1])
\begin{align}\label{eq:Th_adjoint-u-_system}
    \begin{split}
\dot{\eta}^{\bar{u}}(t) &= - \left( \frac{\partial f}{\partial x} (x^{(\bar{u},\bar{d})}(t), \bar{u}(t), \bar{d}(t) ) \right)^\top \eta^{\bar{u}}(t), \\  \eta^{\bar{u}}(\bar{t}) &= -\left(D g_{i^*}(z) \right)^\top
    \end{split}
    \end{align}
    where the orientation of the normal $\left(D g_{i^*}(z) \right)$ to $\partial X_{\bar{t}}(\bar{x},\bar{u})$ can be chosen arbitrarily. Moreover,
    \begin{align}\label{eq:Th_Hamilton-u-_system}
    \begin{split}
\min_{d\in \DD}  \eta^{\bar{u}}(t)^\top f(x^{(\bar{u},\bar{d})}(t), \bar{u}(t), d ) = \eta^{\bar{u}}(t)^\top f(x^{(\bar{u},\bar{d})}(t), \bar{u}(t), \bar{d}(t) ) =0
 \end{split}
    \end{align}
for almost all $t\leq \bar{t}$, 
where $\bar{t}$ denotes the time at which $z = x^{(\bar{u},\bar{d})}(\bar{t})$ is reached, with $z \in \G_0$ satisfying
    \begin{align}\label{eq:Th_UT-u-Condition}
    \begin{split}
    g_i(z) &= 0, \; \; i \in \II(z), \\
    \min_{d \in \DD}  \max_{i \in \II(z)} L_f g_i(z,\bar{u},d) 
	&= L_f g_{i^*}(z,\bar{u},\bar{d}) = 0.
     \end{split}
    \end{align}
Then, using the fact that 
	\begin{align*}
		\min_{d\in \DD}  \eta^{\bar{u}}(t)^\top f(x^{(\bar{u},\bar{d})}(t), \bar{u}(t), d ) = \max_{d\in \DD}  \eta^{\bar{u}}(t)^\top(- f(x^{(\bar{u},\bar{d})}(t), \bar{u}(t), d )) = \max_{d\in \DD} - \eta^{\bar{u}}(t)^\top f(x^{(\bar{u},\bar{d})}(t), \bar{u}(t), d )
	\end{align*}
and noting that if $- \eta^{\bar{u}}$ is a solution of \eqref{eq:Th_adjoint-u-_system}, then $\eta^{\bar{u}}$ is also a solution of the same equation, we immediately see that $- \eta^{\bar{u}}= \eta^{\bar{d}}= \lambda^{(\bar{u},\bar{d})}$, since they satisfy the same final condition.
The equivalence of \eqref{eq:Th_Hamilton-u-_system}-\eqref{eq:Th_Hamilton-d-_system} with \eqref{eq:Th_Hamiltonian} readily follows.
\end{proof}

\begin{remark}\label{rem:algo}
Based on Theorem~\ref{Th:Barrier_Theorem}, the barrier can be computed as follows.
\begin{enumerate}
\item Candidate points on the active constraint set $\G_0$ are determined by solving  $\max_{i= 1, \ldots, p} g_{i}(z) =0$  and picking the $i's$ for which this maximum is achieved. The corresponding set of such $z$'s is locally a $C^1$ $(n-1)$-dimensional manifold.
\item Among the $z$'s, select those, as well as $u$ and $d$, such that $f$ is orthogonal to $Dg_{i}(z)$ for any $i\in \II(z)$ according to \eqref{eq:Th_UT-Condition} (recall from Remark~\ref{rem:H6} that 
if $\II(z)$ contains more than 1 element, $\max_{i \in \II(z)} g_{i}(z) =0$ implies that all $Dg_{i}(z)$ for $i\in \II(z)$ must be equal). Thus, $z$ now belongs to a locally $(n-2)$-dimensional manifold, that should be $C^1$-diffeomorphic to a polytope of $\R^{n-2}$ to satisfy (H6).
\item Next, the minimizing control and maximizing disturbance laws are obtained from \eqref{eq:Th_Hamiltonian}, as functions of time, the final $\bar{t}$ being arbitrary.
\item Finally, the system dynamics, with $u$ and $d$ obtained at the previous step, are integrated backwards in time from $\bar{t}$ and every $z$ obtained above, together with the adjoint equation \eqref{eq:Th_adjoint_system}, resulting in trajectories that are candidate to run in the barrier. From this integration, the choice of orientation of the normal $Dg_{i}(z)$ may be possibly made in order to guarantee that the resulting integral curves well belong to $\G_-$.
\end{enumerate}
We stress that the obtained integral curves are only \emph{candidate} barrier trajectories satisfying necessary conditions, not sufficient in general,  and a final verification of their optimality has to be made by inspection.
\end{remark}

\begin{remark}
While HJB and reachability methods compute trajectories and extremal feedback controllers, interior to the $n$-dimensional admissible set, the barrier method focuses exclusively on the $(n - 1)$-dimensional boundary. This way, it avoids solving PDEs whose natural complexity is even increased by the state constraints while the barrier is computed through a set of ordinary differential equation where the constraints are absent.
\end{remark}
\section{Back to the Example of Adaptive Cruise Control}\label{Sec:Example}

We now return to the adaptive cruise control example introduced in Section~\ref{Sec:Introduction}.
We recall that the system dynamics describing two vehicles driving in convoy are given by \eqref{eq:ACC}, together with the state constraints \eqref{eq:ACC_constraints}, which enforce a minimum and maximum distance between the vehicles.

We briefly verify that the assumptions (H1)-(H4) are satisfied for the considered system.
Assumption (H1) holds since the system dynamics are polynomial in the state and thus $\mathcal{C}^{\infty}$.
Assumption (H4) is satisfied because the constraint functions $g_1$ and $g_2$ are affine, and their respective zero level sets define smooth $(n-1)$-dimensional manifolds.
Assumption (H3) holds since the system is affine in the control and disturbance inputs, and the sets $\UU$ and $\DD$ are convex. Assumption (H2), however, is not satisfied globally due to the quadratic term in the dynamics of $\dot{x}_2$.
Nevertheless, the analysis is restricted to a compact subset of the state space, so that the required boundedness properties of trajectories remain valid.
This restriction is natural in the present setting, as the state constraints impose bounds on the inter-vehicle distance $x_3$, and the combination of drag effects in the dynamics and bounded inputs prevents from an unbounded growth of the velocities $x_1$ and $x_2$, so that the relevant trajectories evolve in a compact region of the state space.

\subsection{Constraints dominated by $g_1$}\label{SubSec:g1_dominating}
We first analyze the barrier trajectories when $g_1(x) \geq g_2(x)$ in a neighborhood of a point $z \triangleq (z_1,z_2,z_3)$ with $g_1(z) = 0$.
The ultimate tangentiality condition \eqref{eq:Th_UT-Condition} at this point $z$ reads
\begin{align*}
   \min_{u \in \UU} \, \max_{d \in \DD} \, L_f g_1(z,u,d) &= \min_{u \in \UU} \, \max_{d \in \DD} \, \begin{pmatrix}
        0 & \tau & -1
    \end{pmatrix} \begin{pmatrix}
        a + d_1 \\
        -\left(a_0 + a_1 z_2 + a_2 z_2^2\right) + g d_2 + g u \\
        z_1 - z_2
    \end{pmatrix} \\
    &= \min_{u \in \UU} \, \max_{d \in \DD} \, -\tau \left( a_0 +  a_1 z_2 + a_2 z_2^2 - g d_2 - g u\right) - z_1 + z_2 \\
    &= -\tau a_0 - \tau a_1 z_2 - \tau a_2 z_2^2 + \tau g \overline{d}_2 + \tau g \underline{u} - z_1 + z_2 \\
    &= - \tau a_2 z_2^2 + \left(1 - \tau a_1 \right) z_2 - z_1 - \tau a_0 + \tau g \underline{u} + \tau g \overline{d}_2 \\
    & = 0.
\end{align*}
Thus, assuming that $(1 -\tau a_1)^2 - 4 \tau a_2 (z_1 + \tau a_0 - \tau g \underline{u} - \tau g \overline{d}_2) \geq 0$, we have
\begin{align*}
    0 &= z_2^2 + \frac{a_1 \tau - 1}{\tau a_2} z_2 + \frac{z_1 + \tau (a_0 - g \underline{u} - g \overline{d}_2)}{\tau a_2} \\
    \Leftrightarrow z_2^{\pm} &= \frac{1 -\tau a_1}{2 \tau a_2} \pm \sqrt{\left(\frac{1 -\tau a_1}{2 \tau a_2}\right)^2 \! \! - \frac{z_1 + \tau (a_0 - g \underline{u} - g \overline{d}_2)}{\tau a_2}} .
\end{align*}
If the product of the roots is non positive, \ie, $z_1 + \tau a_0 - \tau g \underline{u} - \tau g \overline{d}_2 \leq 0$, we obtain two solutions with opposite signs and, from the condition $g_1(z) = 0$, we get $z_3 = \tau z_2$.

We will see later on, when simulating the barrier trajectories, that the ones starting from the positive root must be discarded since they approach the constraints from the complement of the constrained state space and hence are not part of the robust admissible set.

The remaining part of the barrier must satisfy \eqref{eq:Th_Hamiltonian}:
\begin{align*}
     0 & = \min_{u \in \UU} \, \max_{d \in \DD} \, (  \lambda(t))^\top f(x(t), u,d) \\
     &= \min_{u \in \UU} \, \max_{d \in \DD} \, (\lambda(t))^\top \begin{pmatrix}
        a + d_1 \\
        -\left(a_0 + a_1 x_2(t) + a_2 x_2^2(t)\right) + g d_2 + g u \\
        x_1(t) - x_2(t)
    \end{pmatrix} \\
    &= \min_{u \in \UU} \, \max_{d \in \DD} \, \lambda_1(t) (a + d_1) + \lambda_3(t) (x_1(t) - x_2(t)) + \lambda_2(t) \left( -a_0 - a_1 x_2(t) - a_2 x_2^2(t) + g d_2 + g u \right) ,
\end{align*}
from which we can deduce the barrier control law 
\begin{align}\begin{split}\label{eq:ubar}
    u(t) = \begin{cases}
        \overline{u} & : \lambda_2(t) < 0 \\
        \underline{u} & : \lambda_2(t) > 0 \\
        0 & : \lambda_2(t) = 0
    \end{cases}
    \end{split}
\end{align}
and the optimal disturbances
\begin{align}\begin{split}\label{eq:dbar}
    d_1(t) = \begin{cases}
        \overline{d}_1 & : \lambda_1(t) > 0 \\
        \underline{d}_1 & : \lambda_1(t) < 0 \\
        0 & : \lambda_1(t) = 0
    \end{cases}, \; \; d_2(t) = \begin{cases}
        \overline{d}_2 & : \lambda_2(t) > 0 \\
        \underline{d}_2 & : \lambda_2(t) < 0 \\
        0 & : \lambda_2(t) = 0
    \end{cases}.
    \end{split}
\end{align}
Here, we use the convention that for $\lambda_i(t) = 0$, $i = 1,2$, the respective control and disturbance values are chosen as zero, although they can have any (admissible) value in this case.
The adjoint $\lambda$ evolves over time according to the adjoint equation \eqref{eq:Th_adjoint_system}, which reads
\begin{align}\label{eq:ACC_adjoint}
\begin{split}
    \dot{\lambda}(t) &= - \begin{pmatrix}
        0 & 0 & 0 \\
        0 & -2 a_2 x_2(t) - a_1 & 0 \\
        1 & -1 & 0
    \end{pmatrix}^\top \lambda(t) = \begin{pmatrix}
        - \lambda_3(t) \\
        \left( 2 a_2 x_2(t) + a_1 \right) \lambda_2(t) + \lambda_3(t) \\
        0
    \end{pmatrix}. \\
    \lambda(0) &= (D g(z))^\top = \begin{pmatrix}
        0 & \tau & -1
    \end{pmatrix}^\top .
\end{split}
\end{align}
Note that, as remarked before, taking into account their final condition, the barrier trajectories are indeed parameterized by $z_1$, the final speed of the leading vehicle.

\subsection{Constraints dominated by $g_2$}\label{SubSec:g2_dominating}
We now consider the case $g_1(x) \leq g_2(x)$ in a neighborhood of $z$ with $g_2(z) = 0$.
The ultimate tangentiality condition \eqref{eq:Th_UT-Condition} reads:
\begin{align*}
    \min_{u \in \UU} \, \max_{d \in \DD} \, L_f g_2(z,u,d) = \min_{u \in \UU} \, \max_{d \in \DD} \, \begin{pmatrix}
        0 & 0 & 1
    \end{pmatrix} \! \begin{pmatrix}
        a + d_1 \\
        -\left(a_0 + a_1 z_2 + a_2 z_2^2 \right) + g d_2 + g u \\
        z_1 - z_2
    \end{pmatrix} = \min_{u \in \UU} \, \max_{d \in \DD} \, z_1 - z_2 = 0,
\end{align*}
thus $z_2 = z_1$, and with $g_2(z) = 0$, we obtain $z_3 = D_{\textnormal{max}}$.
Thus, the ultimate tangentiality points in this case are $(z_1, z_1, D_{\textnormal{max}})^\top$, again parameterized by $z_1$.

The minimizing control, maximizing disturbances and the adjoint system equations are the same as in the first case with dominating constraint $g_1$.
Only the final value of the adjoint\footnote{As previously mentioned, for simplicity's sake the final time is set equal to $0$.} differs and is $\lambda(0) = (Dg_2(z))^\top = (0,0,1)^\top$.

\subsection{Reparameterization with Respect to the Leading Vehicle's Speed}\label{SubSec:Reparameterization}
Before simulating the barrier trajectories of system \eqref{eq:ACC}, we need to determine the physically reasonable domain on which we want to perform the analysis of the system, namely the range of the uncontrolled state $x_1$, that we assume $x_1 \geq 0$.

We do not allow the first car to accelerate indefinitely, since the second car's acceleration is limited by the air resistance at some point.
Thus, if car 1 accelerates constantly, the upper constraint $g_2$ on the distance between the vehicles will always be violated. 
We therefore assume that $x_1\leq \hat{x}_1$,
this upper bound $\hat{x}_1$ being supposed to be compatible with the second car's acceleration.
According to \eqref{eq:ACC}, the largest admissible acceleration of the first car conditionally to the largest worst-case acceleration of the second one is given by
\begin{align*}
    0 &= \dot{x}_1 - \dot{x}_2
   = a + \overline{d}_1 + \left(a_0 + a_1 \hat{x}_1 + a_2 \hat{x}_1^2\right) - g \underline{d}_2 - g \overline{u} \\
    \Leftrightarrow \hat{x}_1^{\pm} &= \frac{1}{2 a_2} \left( - a_1 \pm \sqrt{a_1^2 + 4 a_2 (g d_2 + g u - a_0 - a - \overline{d}_1)} \right),
\end{align*}
where the negative root will be discarded later.
The remaining one is denoted by $\hat{x}_1^{+}$. Therefore we consider that the system state $(x_{2},x_{3})$ is parameterized by $x_{1}\in [0, \hat{x}_1^{+}]$.

This suggests that we re-parameterize the whole system dynamics by $x_1$  in place of the time, \ie
$$\frac{d}{dx_{1}}= \frac{1}{\frac{dx_1 }{d t}} \frac{d}{dt}= \frac{1}{a+d_{1}}\frac{d}{dt}$$
Thus
 \begin{align}\begin{split}\label{eq:z1paramsys}
 \frac{d x_2}{d x_1} &= \frac{-\left(a_0 + a_1 x_2 + a_2 x_2^2\right) + g d_2 + g u}{a + d_1}\\
\frac{d x_3}{d x_1} &= \frac{x_1 - x_2}{a + d_1}.
\end{split}
\end{align}

The reader may verify that the reparameterized barrier control \eqref{eq:ubar} and disturbance \eqref{eq:dbar} remain unchanged, while the new adjoint system reads
\begin{align}\begin{split}\label{eq:z1paramad}
    \frac{d \mu}{d x_1} = \begin{pmatrix}
        \frac{2 a_2 x_2 + a_1}{a + d_1} & \frac{1}{a + d_1}\\
        0 & 0
    \end{pmatrix} \begin{pmatrix}
        \mu_1 \\ \mu_2
    \end{pmatrix}.
    \end{split}
\end{align}

Forward integration with respect to $t$ means that $x_1$ increases if $a+d_1 \geq 0$ and decreases otherwise.
Accordingly, backward integration with respect to $t$ means that $x_1$ decreases if $a+d_1 \geq 0$ and increases otherwise. 
Thus, in the case where the constraint is dominated by $g_{1}$, we integrate the system \eqref{eq:z1paramsys} and its adjoint system \eqref{eq:z1paramad} over the interval $[z_1, \hat{x}_1^{+}]$ for each ultimate tangentiality point $z$ parameterized by $z_1 \in [0,\hat{x}_1^{+}]$.
In the case where the constraint is dominated by $g_2(z)$, the integration is done over the interval $[z_1,0]$, also for each ultimate tangentiality point $z$ parameterized by $z_1 \in [0,\hat{x}_1^{+}]$.

Should the sign of $\dot{x}_1$ on the barrier trajectories change over time, this integration has to be done separately  on each interval where the sign of $\dot{x}_1$ is constant. Remark that this computation is aimed at providing an off-line (open-loop) construction of the barrier and therefore a real-time dependance on $\dot{x}_1$, including the detection aspects of its sign changes, is not considered here.

\subsection{Simulation Results}\label{SubSec:Simulation_Result}
We use the same parameter configurations as in \cite{Ames_2014} and \cite{Xu_2015} with the constraint parameters $\tau = 1.8$ and $D_{\textnormal{max}} = 100$, as well as the model parameters $m = 1650$, $F_0 = 0.1$, $F_1 = 5$ and $F_2 = 0.25$.
For the input constraints, we choose $\UU = [-0.5,0.5]$, $\DD_1 = [-0.3, 0.3]$ and $\DD_2 = [-0.4,0.4]$.
Moreover, we assume that the leading vehicle aims at driving at a constant speed, i.e., $a = 0$.

\begin{figure}[htbp]
\centerline{\includegraphics[width=90mm,draft=false]{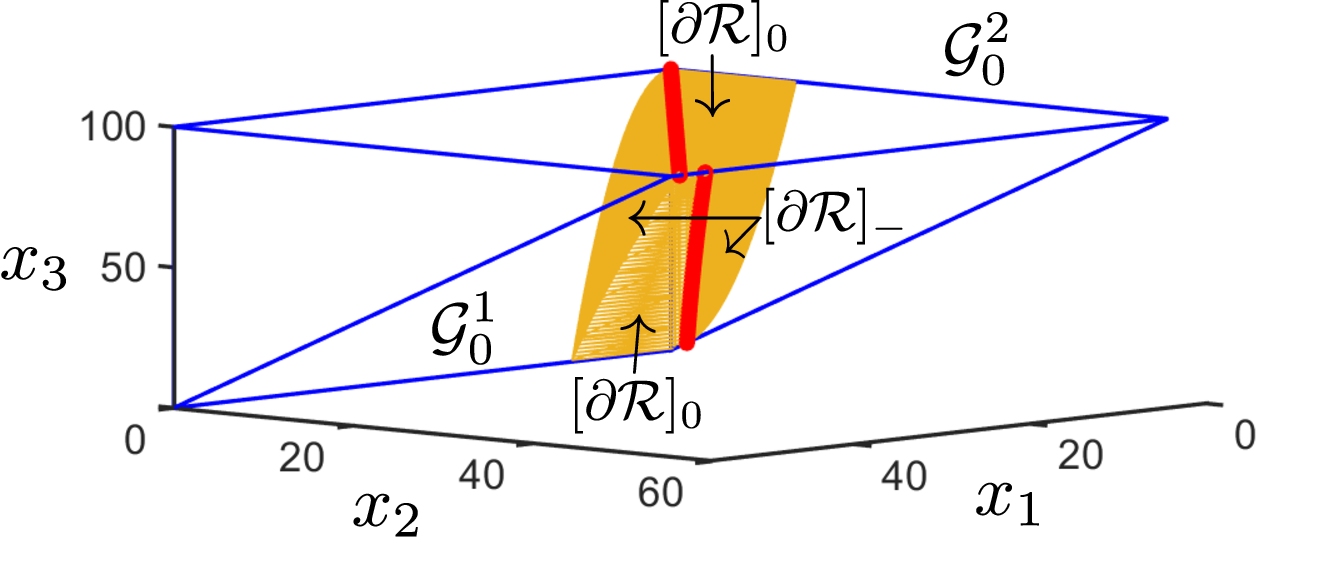}}
\caption{The robust admissible set of system \eqref{eq:ACC} in the $3$-dimensional $x_1$-$x_2$-$x_3$-space in orange with points of ultimate tangentiality in red and state constraints \eqref{eq:ACC_constraints} outlined in blue.\label{fig:ACC_3d-Plot}}
\end{figure}

\begin{figure}[htbp]
\centerline{\includegraphics[width=78mm,draft=false]{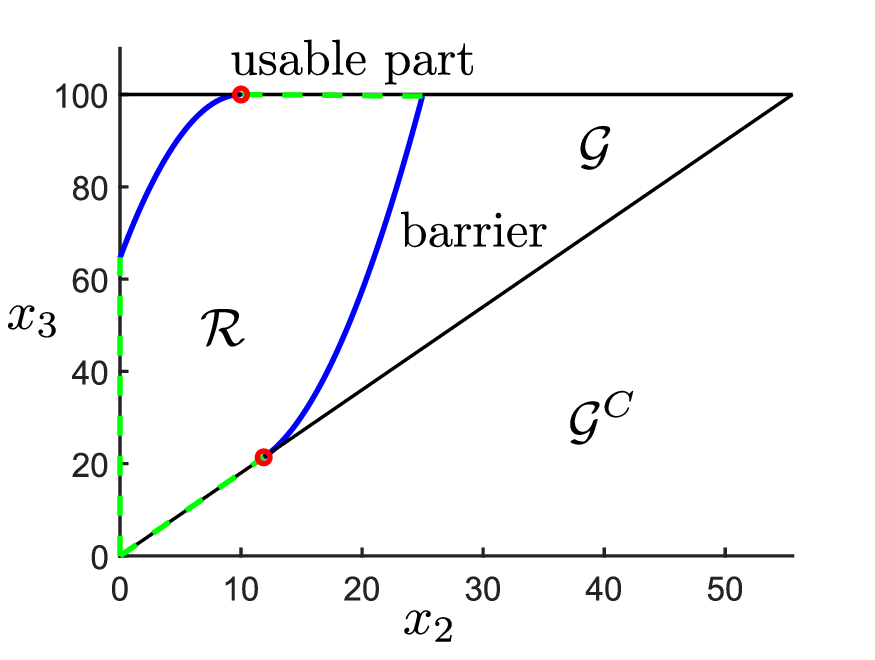}}
\caption{Barrier trajectories of system \eqref{eq:ACC} parameterized for $z_1 = 10$ in the $x_2$-$x_3$-plane in blue with points of ultimate tangentiality in red and usable part in green. State constraints \eqref{eq:ACC_constraints} in black.\label{fig:ACC_x2-x3-Plot_z1=10}}
\end{figure}

\begin{figure}[htbp]
\centerline{\includegraphics[width=78mm,draft=false]{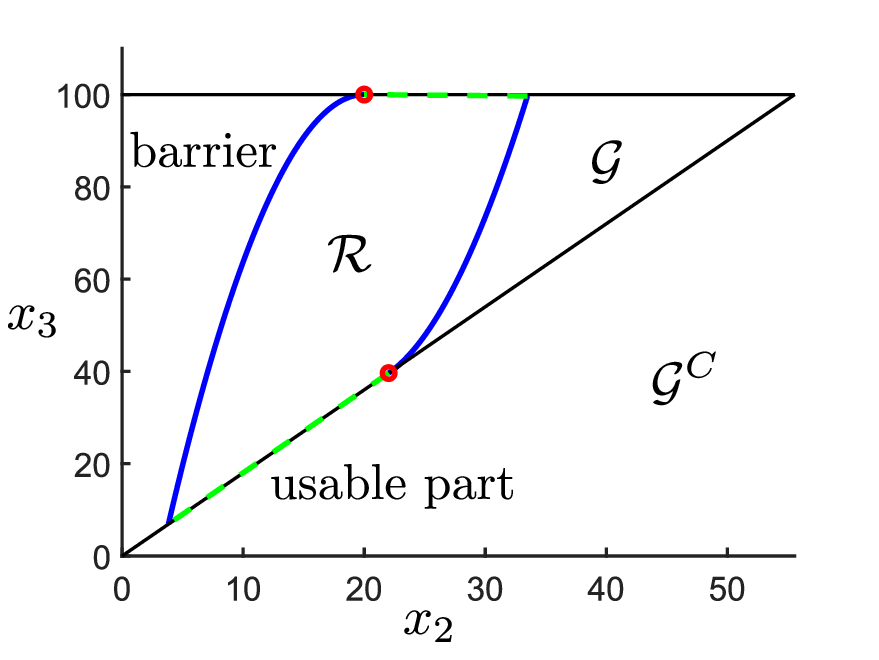}}
\caption{Barrier trajectories of system \eqref{eq:ACC} parameterized for $z_1 = 20$ in the $x_2$-$x_3$-plane in blue with points of ultimate tangentiality sampled in red and usable part in green. State constraints \eqref{eq:ACC_constraints} in black.\label{fig:ACC_x2-x3-Plot_z1=20}}
\end{figure}

\begin{figure}[htbp]
\centerline{\includegraphics[width=78mm,draft=false]{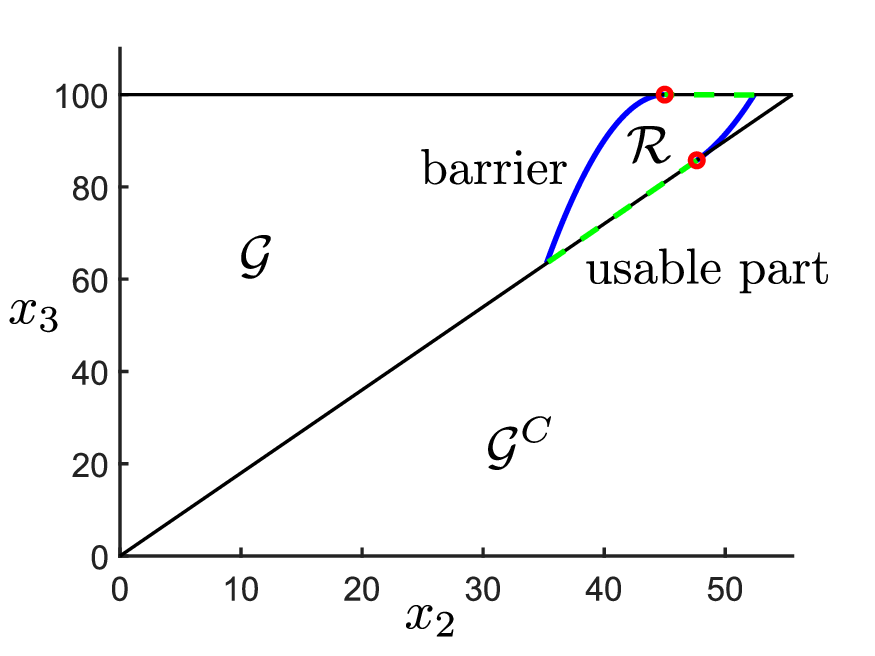}}
\caption{Barrier trajectories of system \eqref{eq:ACC} parameterized for $z_1 = 45$ in the $x_2$-$x_3$-plane in blue with points of ultimate tangentiality sampled in red and usable part in green. State constraints \eqref{eq:ACC_constraints} in black.\label{fig:ACC_x2-x3-Plot_z1=45}}
\end{figure}

The resulting robust admissible set is shown in Figure~\ref{fig:ACC_3d-Plot} and projections for specific values of $z_1$ are shown in Figures~\ref{fig:ACC_x2-x3-Plot_z1=10}-\ref{fig:ACC_x2-x3-Plot_z1=45} for the two constraints $g_{1}$ and $g_{2}$ simultaneously, with the points of ultimate tangentiality in red, the barrier trajectories in blue, the usable part in green and the state constraints in black,  for $z_1 = 10$, $z_1 = 20$ and $z_1 = 45$.

The robust admissible set $\RA$ is the area delimited by the barrier trajectories and the state constraints.
If the state is located in $\RAC\cap \G$, though currently still safe, the system is guaranteed to  violate the state constraints in finite time in the future under worst-case disturbances.
Note that if the disturbances are not worst-case, the system might still be able to enter back into $\RA$, although this cannot be guaranteed.
Intuitively speaking, the system becomes all the "more unsafe" the further it is away from $\RA$ and the closer to $\G_0$.

We may also remark that the two cases of dominating constraints are here disjoint. Nevertheless, in general, nothing prevents from these cases to overlap. In such overlapping points, since the constraints must be equal, they locally also have equal gradients since they realize the maximum of the constraints. Hence, the ultimate tangentiality conditions corresponding to the dominating constraints are equal.

It can be observed that for higher speeds of the first vehicle at the ultimate tangentiality points, the robust admissible set is smaller.
This effect occurs partially due to the state constraint $g_1$ requiring a larger physical distance between the vehicles while the upper constraint given by $g_2$ remains constant.
This can be observed, for instance, when comparing the size of the usable part in figures \ref{fig:ACC_x2-x3-Plot_z1=10} and \ref{fig:ACC_x2-x3-Plot_z1=45}.

Note that in  \cite{Ames_2014,Xu_2015}, a similar example is presented with the only constraint $g_1$ and with disturbances on the follower only. 
The approach of these authors consisted in computing a control barrier function or, in other words, to design a feedback ensuring the asymptotic convergence of the follower to a subset $\cal{C}\subset\RA$, interpreted as the desired distance to the leader with a desired speed $v_{d}$, given the disturbances. Therefore the control barrier function framework appears to be complementary to ours, providing information of a different nature: in their approach no precise indication of the width of the admissible set is obtained whereas our approach does not provide a feedback to guarantee by itself the constraint satisfaction.

\section{Conclusions}\label{Sec:Conclusion}

In this paper, we extended the results related to the boundary of the admissible set from \cite{Levine_2013} to include disturbances and the maximal robust positively invariant set from \cite{Esterhuizen_2020} to include controls. In the setting of robust admissible sets, both control and disturbance inputs act simultaneously on the state of a nonlinear system subject to state and input constraints. We have obtained a characterization of the boundary of the robust admissible set, divided in two parts, the usable part and the barrier. We have proven that the barrier is made of system trajectories running along it until they intersect tangentially with the active state constraints. Moreover, the corresponding tangent hyperplane is locally separating the robust admissible set and its complement, and is given by a saddle-point condition of the Hamiltonian along the barrier, a saddle-point version of the Pontryagin's maximum principle.
These theoretical results were then applied to compute the boundary of the robust admissible set for an adaptive cruise control system, originally introduced in \cite{Ames_2014} and extended to a robust setting in \cite{Xu_2015}. The control barrier function framework adopted by these authors appears to be complementary to ours.

Let us stress, moreover, that the example shows that our barrier construction does not need to be re-computed for each parameter value since it is done once for all exploiting an explicit dependence on the system dynamics.

Our approach might open the possibility of solving optimal control or model predictive control problems subject to state and input constraints and to disturbances in two steps by performing the optimization on the state space restricted to the robust admissible set, once its boundary is explicitly known. We thus may evaluate if some simplifications can be obtained in this way. Note that it requires to have access to a suitably simple numerical representation of the admissible set to avoid deteriorating the performances of such an approach. Future research may address this challenge.

\section*{Acknowledgments}

The authors would like to thank Dr.~Willem Esterhuizen for fruitful discussions on this topic.

This work was funded by the Deutsche Forschungsgemeinschaft (DFG, German Research Foundation) – Project Number 531896505.

\bibliographystyle{plain}
\bibliography{wileyNJD-AMS}

\appendix

\section{Needle Perturbations}\label{Appendix:Needle_Perturbations}
Consider a control $\bar{u} \in \U$ and a disturbance $\bar{d} \in \D$, as well as the combined input $\bar{v} = (\bar{u},\bar{d}) \in \U \times \D$.
Moreover, for some bounded $\varepsilon_0 > 0$, let $\varepsilon \in [0,\varepsilon_0]$ and consider two constants $T,L \in \R$.
For $\kappa \triangleq (\tilde{v},\tau,l) \in (\UU \times \DD) \times [0,T] \times [0,L]$ with $\tilde{v} = (\tilde{u},\tilde{d})$, a variation $v_{\kappa,\varepsilon}$ of the combined input $\bar{v}$ is given by
\begin{align}\label{eq:Def_needle_perturbation}
    v_{\kappa,\varepsilon} \triangleq \bar{v} \Join_{(\tau - l \varepsilon)} \tilde{v} \Join_\tau \bar{v} = \begin{cases}
        \tilde{v} & \mathrm{on~} [\tau - l \varepsilon, \tau), \\
        \bar{v} & \mathrm{elsewhere~on~} [0,T] .
    \end{cases}
\end{align}
Note that $\tilde{v} \in \UU \times \DD$ is a constant input over the interval $[\tau - l \varepsilon)$.

Additionally, consider an initial state perturbation $h \in \R^n$ with $\| h \| \leq H$ for some $H \in \R$.
Thus, the corresponding integral curve $x^{(\bar{x}+ \varepsilon h, v_{\kappa,\varepsilon})}$ is starting from $\bar{x}+ \varepsilon h$ with the perturbed combined input $v_{\kappa,\varepsilon}$.

Note that for all $t \in [0, \tau - l \varepsilon)$ we have $x^{(\bar{x}+ \varepsilon h, v_{\kappa,\varepsilon})}(t) = x^{(\bar{x}+ \varepsilon h, \bar{v})}(t) = x^{(\bar{x}+ \varepsilon h, \bar{u}, \bar{d})}(t)$.
Furthermore, by denoting $z_\varepsilon(\tau - l \varepsilon) \triangleq x^{(\bar{x} + \varepsilon h, \bar{v})}(\tau - l \varepsilon)$ and $z_\varepsilon(\tau) \triangleq x^{(\bar{x} + \varepsilon h, v_{\kappa, \varepsilon})}(\tau)$, we also have
\begin{align*}
    x^{(\bar{x} + \varepsilon h, 0, v_{\kappa, \varepsilon})} = x^{(\bar{x} + \varepsilon h, 0, \bar{v})} \Join_{\tau - l \varepsilon} x^{(z_\varepsilon(\tau - l \varepsilon), \tau - l \varepsilon, \tilde{v})} \Join_\tau x^{(z_\varepsilon(\tau), \tau, \bar{v})}.
\end{align*}
Lastly, the fundamental matrix of the variational equation is given by 
\begin{align*}
    \frac{d}{dt} \Phi^{(\bar{u},\bar{d})}(t,s) &= \left( \frac{\partial f}{\partial x} (x^{(\bar{x},\bar{u},\bar{d})}(t),\bar{u}(t),\bar{d}(t)) \right) \Phi^{(\bar{u},\bar{d})}(t,s), \\
    \Phi^{(\bar{u},\bar{d})}(s,s) &= I_n,
\end{align*}
where $I_n$ denotes the identity matrix in $\R^{n \times n}$.
Similar to \cite{Levine_2013}[Lemma B.1], we can now formulate the following approximation result from, \eg, \cite{Pontryagin_1962}[Chapter II, Section 13], \cite{Lee_1967}[Chapter 4, page 248] or \cite{Gamkrelidze_1999}, which is adjusted to account for the combined input $\bar{v}$.

\begin{lemma}\label{Lem:Needle_Perturbation_Approx}
    The sequence $\big\lbrace x^{(\bar{x}+\varepsilon h, v_{\kappa,\varepsilon})} \big\rbrace_{\varepsilon \geq 0}$ is uniformly converging to $x^{(\bar{x},\bar{u},\bar{d})}$ on $[0,T]$ as $\varepsilon$ tends to $0$, uniformly with respect to $h$ and $\kappa$.

    Moreover, if $\tau$ is a Lebesgue point of $\bar{u}$ and $\bar{d}$, we have for all $t \in [\tau, T]$
\begin{align}\label{eq:Lemma_Needle_Perturb_Approx}
\begin{split}
     x^{(\bar{x}+\varepsilon h, v_{\kappa,\varepsilon})}(t) - x^{(\bar{x},\bar{v})}(t) = x^{(\bar{x}+\varepsilon h, v_{\kappa,\varepsilon})}(t) - x^{(\bar{x},\bar{u},\bar{d})}(t) = \varepsilon w(t,\kappa,h) + O(\varepsilon^2),
 \end{split}
\end{align}
where
\begin{align}\label{eq:Lemma_Needle_Perturb_w}
\begin{split}
	w(t,\kappa,h) &\triangleq \Phi^{(\bar{u},\bar{d})}(t,0) h + l \Phi^{(\bar{u},\bar{d})}(t,\tau) \, \cdot \left( f(x^{(\bar{x},\bar{v})}(\tau),\tilde{v}) - f(x^{(\bar{x},\bar{v})}(\tau),\bar{v}(\tau)) \right) \\
	&= \Phi^{(\bar{u},\bar{d})}(t,0) h + l \Phi^{(\bar{u},\bar{d})}(t,\tau) \, \cdot \left( f(x^{(\bar{x},\bar{u},\bar{d})}(\tau),\tilde{u},\tilde{d}) - f(x^{(\bar{x},\bar{u},\bar{d})}(\tau),\bar{u}(\tau),\bar{d}(\tau)) \right) .
     \end{split}
\end{align}
\end{lemma}
Note that the perturbed vector $w(t,\kappa,h)$ is tangent at $x^{(\bar{x},\bar{v})}(t)$ to the curve $\varepsilon \longmapsto x^{(\bar{x} + \varepsilon h, v_{\kappa, \varepsilon})}(t)$.

Moreover, by multiplying $l$ or $h$ with any real number $\mu > 0$, we obtain the perturbation vector $\mu \cdot w(t,\kappa,h)$.
Thus, the sets of the so-obtained perturbation vectors form a cone in $\R^n$.

\section{The Perturbation Cone}\label{Appendix:Perturbation_Cone}
We extend the variation \eqref{eq:Def_needle_perturbation} to a sequence of variations, given by the vector $\chi \triangleq \lbrace (\tilde{u}_i,\tau_i,l_i) \, \vert \, i = 1,2, \ldots k \rbrace$ such that
\begin{align*}
    v_{\chi,\varepsilon} \triangleq \bar{v} \Join_{(\tau_1 - l_1 \varepsilon)} \tilde{v}_1 \Join_{\tau_1} \bar{v} \Join_{(\tau_2 - l_2 \varepsilon)} \tilde{v}_2 \Join_{\tau_2} \bar{v} \ldots \Join_{(\tau_k - l_k \varepsilon)} \tilde{v}_k \Join_{\tau_k} \bar{v},
\end{align*}
where the perturbation times $\tau_1 < \tau_2 < \ldots < \tau_k$ are Lebesgue points of both, $\bar{u}$ and $\bar{d}$.
The statement of Lemma \ref{Lem:Needle_Perturbation_Approx} still holds, and we obtain for every $t \in [\tau_k,T]$ the approximation \eqref{eq:Lemma_Needle_Perturb_Approx} as
\begin{align*}
    x^{(\bar{x} + \varepsilon h,v_{\chi,\varepsilon})}(t) - x^{(\bar{x},\bar{v})}(t) = \varepsilon w(t,\chi,h) + O(\varepsilon^2)
\end{align*}
with
\begin{align*}
    w(t,\chi,h) \triangleq \Phi^{(\bar{u},\bar{d})}(t,0) h + \sum_{i = 1}^k l_i \Phi^{(\bar{u},\bar{d})}(t,\tau_i) \, \cdot  \left(f(x^{(\bar{x},\bar{v})}(\tau_i),\tilde{v}_i) - f(x^{(\bar{x},\bar{v})}(\tau_i),\bar{v}(\tau_i)) \right).
\end{align*}
The variation $v_{\chi,\varepsilon}$ of the combined input $\bar{v}$ is a so-called \textit{needle perturbation}.
We denote the cone generated by convex combinations of these needle perturbations at a time $t \in [0,T]$ by $\mathcal{K}_t$.

We can perturb the control $\bar{u}$ and the disturbance $\bar{d}$ independently by choosing $\tilde{d}_i = \bar{d}(\tau_i)$ and $\tilde{u}_i = \bar{u}(\tau_i)$ for all $i = 1,2,\ldots, k$, respectively.
In order to easily distinguish these specified cones, we denote the (convex) cone of needle perturbations of the control input $\bar{u}$ at time $t \in [0,T]$ by $\mathcal{K}^{\bar{u}}_t$ and the cone related to a perturbation of the disturbance input $\bar{d}$ by $\mathcal{K}^{\bar{d}}_t$.

Note that for all $0 \leq \tau_1 \leq \tau_2 \leq T$ we have $\Phi^{(\bar{u},\bar{d})}(\tau_{2},\tau_{1}) \mathcal{K}_{\tau_1} \subset \mathcal{K}_{\tau_2}$, which also directly applies to the cones $\mathcal{K}^{\bar{u}}_t$ and $\mathcal{K}^{\bar{d}}_t$.

\section{The Maximum Principle}\label{Appendix:Maximum_Principle}
\vspace*{12pt}
We follow the approach of \cite{Lee_1967}[Theorem 3, Chapter 4, p. 254] in the proof of Theorem \ref{Th:Barrier_Theorem}.
\begin{theorem}[Maximum Principle]\label{Th:Maximum_Principle}
    Consider system \eqref{eq:control_system} with control input $u \in \U$ and fixed $\bar{d}$.
    Let $\bar{u} \in \U$ be such that $x^{(x_0,\bar{u}, \bar{d})}(t_1) \in \partial X^{\bar{d}}_{t_1}(x_0)$ for some $t_1>0$ where $X^{\bar{d}}_{t}(x_0) \triangleq \{ x^{(x_0,u, \bar{d})}(t) \, \vert \, u \in \U \}$ denotes the reachable set at time $t$ from $x_0$ under the disturbance $\bar{d}$.
    Then, there exists a nonzero absolutely continuous maximal solution $\eta^{\bar{u},\bar{d}}$ to the adjoint equation
    \begin{align}\label{eq:MP_adjoint_system_u}
        \dot{\eta}^{\bar{u},\bar{d}}(t) = - \left( \frac{\partial f}{\partial x} (x^{(x_0,\bar{u}, \bar{d})}(t),\bar{u}(t), \bar{d}(t)) \right)^\top \eta^{\bar{u},\bar{d}}(t)
    \end{align}
    such that
    \begin{align}\label{eq:MP_Hamiltonian_u}
    \begin{split}
		\max_{u \in \UU} \left\lbrace (\eta^{\bar{u}, \bar{d}}(t))^\top f(x^{(x_0,\bar{u}, \bar{d})}(t),u, \bar{d}(t)) \right\rbrace = (\eta^{\bar{u}, \bar{d}}(t))^\top f(x^{(x_0,\bar{u}, \bar{d})}(t),\bar{u}(t), \bar{d}(t)) = \textnormal{constant}
    \end{split}
    \end{align}
    for almost all $t \in [0,t_1]$.
\end{theorem}

We immediately deduce the following:
\begin{corollary}
 Consider system \eqref{eq:control_system} with disturbance $d \in \D$ and fixed $\bar{u}$.
    Let $\bar{d} \in \D$ be such that $x^{(x_0,\bar{u}, \bar{d})}(t_1) \in \partial X^{\bar{u}}_{t_1}(x_0)$ for some $t_1>0$ where $X^{\bar{u}}_{t}(x_0) \triangleq \{ x^{(x_0,\bar{u}, d)}(t) \, \vert \, d \in \D \}$ denotes the reachable set at time $t$ from $x_0$ under the control $\bar{u}$.
    Then, there exists a nonzero absolutely continuous maximal solution $\eta^{\bar{u},\bar{d}}$ to the adjoint equation
    \begin{align}\label{eq:MP_adjoint_system_d}
        \dot{\eta}^{\bar{u},\bar{d}}(t) = - \left( \frac{\partial f}{\partial x} (x^{(x_0,\bar{u}, \bar{d})}(t),\bar{u}(t), \bar{d}(t)) \right)^\top \eta^{\bar{u},\bar{d}}(t)
    \end{align}
    such that
    \begin{align}\label{eq:MP_Hamiltonian_d}
    \begin{split}
		\max_{d \in \DD} \left\lbrace (\eta^{\bar{u}, \bar{d}}(t))^\top f(x^{(x_0,\bar{u}, \bar{d})}(t),u, \bar{d}(t)) \right\rbrace = (\eta^{\bar{u}, \bar{d}}(t))^\top f(x^{(x_0,\bar{u}, \bar{d})}(t),\bar{u}(t), \bar{d}(t)) = \textnormal{constant}
    \end{split}
    \end{align}
    for almost all $t \in [0,t_1]$.

\end{corollary}

\end{document}